\documentclass[ijaa,nonblindrev]{informs4}

\usepackage{eqndefns-left} 
\RequirePackage{tgtermes}
\RequirePackage{newtxtext}
\RequirePackage{newtxmath}
\RequirePackage{bm}
\RequirePackage{endnotes}

\usepackage{caption}
\usepackage{subcaption}
\usepackage{adjustbox}  
\usepackage{hyperref}
\usepackage{booktabs}
\usepackage{float}
\usepackage{algorithm}
\usepackage{algpseudocode}
\usepackage{enumitem}
\usepackage{multirow}
\usepackage{booktabs}
\usepackage{adjustbox}
\usepackage{fix-cm}
\DoubleSpacedXII

\usepackage{algorithm}
\usepackage{algpseudocode}
\usepackage{tikz}

\usepackage{natbib}
 \bibpunct[, ]{(}{)}{,}{a}{}{,}%
 \def\bibfont{\small}%
\usepackage{makecell} 
\hypersetup{
    colorlinks,
    linkcolor={red!50!black},
    citecolor={blue!50!black},
    urlcolor={blue!80!black}
}

\EquationsNumberedThrough    

\TheoremsNumberedThrough     
\ECRepeatTheorems  %

\newcommand{\University}{Cornell University}
\newcommand{\Department}{the ORIE Department}
\newcommand{\university}{cornell}

\newcommand{\revision}[1]{{\color{black}#1}}

\begin{document}


\RUNAUTHOR{Ye et al.}

\RUNTITLE{Cornell Optimizes Final Exam Scheduling}

\TITLE{\University\ Uses Integer Programming to Optimize Final Exam Scheduling}

\ARTICLEAUTHORS{%

\AUTHOR{Tinghan Ye}
\AFF{H. Milton Stewart School of Industrial and Systems Engineering, Georgia Institute of Technology,  \EMAIL{joe.ye@gatech.edu}}
\AFF{School of Operations Research and Information Engineering, Cornell University}

\AUTHOR{Adam S. Jovine}
\AFF{School of Operations Research and Information Engineering, Cornell University, \EMAIL{asj53@cornell.edu}
}
\AUTHOR{Willem van Osselaer}
\AFF{Laboratory for Information and Decision Systems, 
Massachusetts Institute of Technology, \EMAIL{willemvo@mit.edu}}
\AFF{School of Operations Research and Information Engineering, Cornell University}
\AUTHOR{Qihan Zhu}
\AFF{School of Operations Research and Information Engineering, Cornell University,  \EMAIL{qz245@cornell.edu}}
\AUTHOR{David B. Shmoys}
\AFF{School of Operations Research and Information Engineering, Cornell University, \EMAIL{david.shmoys@cornell.edu}}
} 
\ABSTRACT{
This paper presents an integer programming-based optimization framework designed to effectively address the complex final exam scheduling challenges encountered at \University. With high flexibility, the framework is specifically tailored to accommodate a variety of different constraints, including the front-loading of large courses and the exclusion of specific time slots during the exam period. By generating multiple scheduling model variants and incorporating heuristic approaches, our framework enables comprehensive comparisons of different schedules. This empowers the University Registrar to make informed decisions, considering trade-offs in terms of schedule comfort measured by different levels of exam conflicts. Our results demonstrate significant advantage over the historical lecture time-based approach, providing time and effort savings for the university administration while enhancing student and faculty satisfaction.
}%






\maketitle

\section{Introduction} \label{sec:intro}
Academic stress is a universal issue among undergraduate students. And for many, the final exam period is the most stressful time of all. With each exam often requiring hours of studying, it is not only the number of exams but also the timing of exams that decide how challenging and stressful a final exam period will be. Having three exams can be easily manageable if they are spread out over the course of ten days, but three exams on the same day leaves students more stressed and with less time to prepare.

From the perspective of faculty, the best schedule would be one without any scheduling conflicts for students. If a student is scheduled to take two exams in the same time slot, alternative arrangements must be made for one of the exams to be taken at a different time and location. This process can be quite challenging and inconvenient for the faculty. Due to academic integrity concerns, a completely new exam may need to be created specifically for that small group of students.

Because each student takes a different set of classes, finding a final exam schedule means balancing the exam sequences of all the different students; an exam schedule that nicely spreads exams for one student might put all exams in a tight block for another student. When faced with thousands of students taking hundreds of exams, this becomes an impossible task to do by hand.

At \University, faculty and students from \Department\ collaborated with the Office of the University Registrar to schedule final exams. Before using integer programming, final exams at \University\ were scheduled based on a static table that maps first lecture meeting times to exam groups. Courses that started at the same hour would have their finals scheduled  simultaneously. However, the number of meeting times exceeded the available exam slots. Typically, there were around 30 different meeting times but only 22 exam slots, leading to numerous scheduling conflicts. No consideration was given to the other bad events at all. 

During the COVID-19 pandemic, the university sought assistance from \Department\ to address logistics for reopening for residential instruction in Fall 2020 \citep{frazier2022modeling}. Their efforts ranged from completely rescheduling all courses to maximizing in-person attendance while taking into account a 6-foot distance seating requirement \citep{liu2023modeling}, to determining actual classroom seat assignments \citep{greenberg2021automated}. Final exam scheduling was one of the critical tasks undertaken in this collaboration.

We develop an optimization framework based on mixed-integer programming (MIP) to create final exam schedules for the entire \University. The foundation of the framework comprises a multi-stage Group-then-Sequence (GtS) model and an iterative Layer-Cake heuristic algorithm. The framework has been effectively employed for the past five semesters.

\subsection{\revision{Contributions}}

Our proposed methods advance previous research by addressing a range of final exam conflicts. In addition to direct conflicts and back-to-back exams, we consider more complex scenarios, such as having more than two exams within or beyond a 24-hour period. This approach is specifically designed to tackle the unique challenges faced by \University, aiming to enhance both student and faculty comfort during finals week through trade-offs across various conflict levels.

Building upon the existing body of research, this paper presents a novel optimization
framework tailored to the specific needs of students and faculty at \University. The framework offers high flexibility, enabling customization to accommodate constraints such as front-loading large exams or excluding certain time slots. By generating multiple scheduling model variants with different parameter combinations and heuristic approaches, our framework allows for a comprehensive comparison of exam schedules. This comparison helps illustrate trade-offs in schedule comfort, as measured by varying levels of exam conflicts, ultimately assisting the University Registrar in selecting the most optimal final exam schedule.

Our main contribution are as follows:
\begin{enumerate}
    \item We propose a MIP formulation that incorporates penalties for higher-order conflicts. This formulation explicitly addresses scenarios where a student is scheduled for three exams on the same day, within a 24-hour period, or across four consecutive exam slots.
    
    \item We develop an integrated MIP-based framework that synergizes local search and heuristic techniques to efficiently generate exam schedules optimized according to multiple performance metrics. Moreover, our framework supports hybridization strategies to further enhance solution quality.
    
    \item We validate our framework through a comprehensive evaluation at a large American university, where it has been successfully implemented over several consecutive semesters. Additionally, comparative analyses on a benchmark dataset demonstrate that our approach is competitive with existing exam timetabling algorithms,  highlighting its versatility in addressing various exam timetabling challenges.
\end{enumerate}

\subsection{Organization}

The rest of the paper is organized as follows: \revision{First, we detail \University's final exam scheduling problem in Section \ref{sec:prob-descr}. In Section \ref{sec:literature}, we provide a brief review of relevant work.} Section \ref{sec:model} outlines our proposed solution methods, with detailed optimization models and algorithms provided in the Appendix. In Section \ref{sec:implementation}, we present the outcomes of our scheduling models and algorithms using real student enrollment data from \University. Finally, Section \ref{sec:conclsion} concludes our work.

\section{\revision{Problem Description}} \label{sec:prob-descr}

\University\ faces the complex task of scheduling final exams for a significant student body of more than 15,000 and approximately 700 courses at the beginning of every semester. These exams are concentrated within a 7 to 9-day final exam period, with three time slots available per day (9am, 2pm, and 7pm).  However, certain time slots---such as Saturday nights and Sunday mornings---may have restrictions, resulting in roughly 19 to 22 available time slots. The ultimate scheduling problem is to assign around 500 exams to these time slots. The physical locations for the exams are determined separately by the University Registrar, and as such, this study focuses solely on optimizing the timing of the exams. \revision{Room assignment is not a significant concern at \University{} for two primary reasons. First, the university has a sufficient number of classrooms with adequate seating capacities to accommodate demand. The only potential constraint is a cap of 5,000 students taking exams per day; however, all generated schedules remain well below this limit. Second, \University{} typically does not split exams across multiple classrooms, which further simplifies the assignment of physical locations.
}

The goal of an optimized final exam schedule is to enhance the experience for both students and faculty. From the faculty's perspective, student conflicts---such as overlapping or consecutive exams--may require the scheduling of make-up exams according to university policies (\href{https://registrar.\university.edu/calendars-exams/final-exam-policies}{https://registrar.\university.edu/calendars-exams/final-exam-policies}), which imposes additional burdens. Meanwhile, students desire a well-distributed exam schedule to allow for sufficient review and study time. As a result, the optimal exam schedule should minimize various unfavorable metrics, each with its own weight. These metrics include direct conflicts, back-to-back exams, two exams in 24 hours, three exams in a row, and three exams in 4 slots.

 We define the quality of a final exam schedule by the number of ``bad events" it contains:
 \begin{itemize}
     \item Direct conflict: A student with two exams scheduled at the same time slot.
    \item Three exams in a row (triple): A student with three exams assigned to three consecutive slots.
     \item Back-to-back exams (B2B): A student with exams in two consecutive slots.
     \item Two exams in 24 hours (2-in-24hr): A student with two exams in three consecutive slots, or equivalently, within a 24-hour period.
     \item Three exams in 4 slots (3-in-4): A student with three exams in four consecutive slots.
 \end{itemize}
 
We ensure no double counting of metrics; for instance, a student having a pair of exams that contribute to a triple will not be counted again as a B2B or 2-in-24hr.
 
The objective is to minimize the weighted sum of these ``bad events", which serves as a measure of schedule quality. While the primary focus is on eliminating direct conflicts and triples to avoid make-up exams, other metrics are also considered to reduce student stress as much as possible. Additionally, each exam must be assigned to exactly one time slot.

The final exam schedule must also adhere to specific constraints to facilitate room assignments by the University Registrar. In particular, exams for large enrollment courses should be scheduled in earlier time slots to free up larger auditoriums and lecture halls for graduation events and ensure timely grading in line with university deadlines.

\section{\revision{Related Literature}}\label{sec:literature}
This section provides related works in terms of problem formulations and solution techniques.

\subsection{Educational Timetabling}

The scheduling problem encountered by \University\ is situated within the broader domain of educational timetabling.  A comprehensive survey by \cite{ceschia2023educational} offers an up-to-date overview of this evolving field. Educational timetabling is well-supported by an array of benchmark datasets, such as those from the University of Toronto \citep{carter1994general}, the University of Nottingham \citep{burke1996memetic}, and the Purdue University \citep{muller2016real}. These resources, alongside major international competitions like ITC-2007 \citep{mccollum2007second}, ITC-2011 \citep{post2016third}, and ITC-2019 \citep{muller2024real}), as well as tools like UniTime [\href{https://www.unitime.org}{https://www.unitime.org}])), demonstrate the field's practical applications.

Educational timetabling includes three main subfields: high-school timetabling, course timetabling, examination timetabling. Many real-world problems in these areas are approached through integer programming techniques. High-school timetabling, commonly represented in the XHSTT format \citep{post2016third}, has seen the application of MIP in studies like \cite{kristiansen2015integer} and \cite{fonseca2017integer}. Course timetabling problems, such as those at the United States Air Force Academy \citep{gonzalez2018optimal}, the Technion in Israel \citep{strichman2017near}, and MIT \citep{barnhart2022course}, also employ these approaches. \revision{In examination timetabling, there are case studies from the Universidad Politécnica de Madrid \citep{garcia2019universidad} and the 
Department of Education in Vestfold County \citep{avella2022practice}}

\subsection{University Examination Timetabling: Problems and Formulations}

Among educational timetabling problems, our work is most closely related to university examination timetabling problems (ETPs).
Comprehensive surveys of ETPs can be found in the works of \citet{qu2009survey, gashgari2018survey, aldeeb2019comprehensive}, and \citet{siew2024survey}. 

ETPs have been extensively studied for decades, taking various forms with different restrictions based on institutional and regulatory differences. The most relevant formulation for our study is the uncapacitated examination timetabling problem, first introduced by \cite{carter1996examination}. A mathematical model of this formulation is provided by \cite{bellio2021two}. More complex versions that incorporate room assignments and capacity constraints, such as Track 1 of ITC-2007 \citep{mccollum2007second}, exist, but since room capacity is not a factor in our problem, we focus on the uncapacitated version.

In uncapacitated ETPs, the primary hard constraints include avoiding simultaneous exam conflicts and managing the total number of exams assigned to the same time slot \citep{aldeeb2019comprehensive}. Additionally, soft constraints are commonly introduced as penalty terms in the objective function. These typically involve the period spread between  exams with co-enrollments and prioritizing larger exams for earlier scheduling \citep{aldeeb2019comprehensive}. In the literature, many studies focus on metrics related to pairs of exams with co-enrollments, such as back-to-backs, or more broadly speaking, two exams that are $y$ periods apart. Nonetheless, very few of them explicitly address higher-order conflicts, such as scheduling no more than $x$ co-enrolled exams within $y$ consecutive periods ($x, y \geq 3$). This higher-order conflict constraint was first proposed by \cite{carter1996examination} and also appears in the Purdue University formulation \citep{muller2016real}, which includes a soft constraint for a student having three or more exams on the same day, solved using constraint programming. To date, we have not yet encountered MIP formulations that explicitly address such higher-order conflicts.

\subsection{University Examination Timetabling: Solution Techniques}

A comprehensive survey by \cite{aldeeb2019comprehensive} reviews various solution techniques for the uncapacitated ETP. Common approaches include graph-based heuristics, such as those evaluated by \cite{carter1996examination}, and metaheuristics like Simulated Annealing, used in \cite{bellio2021two}.

Our approach to solving this scheduling problem employs MIP for solving subproblems and developing heuristic algorithms. 
Similar mathematical programming-based approaches have been proposed in prior works on general ETPs \citep{woumans2016column, guler2021web, godwin2022obtaining, bazari2023modeling}.

We propose a multi-stage methodology, where a large MIP is divided into smaller stages. Multi-stage approaches have also been explored in the ETP literature \citep{cataldo2017integer, keskin2018examination,al2020practical}. For post-processing, we apply a local search-based algorithm, which shares similarities with existing metaheuristics such as Large Neighborhood Search and fix-and-optimize strategies. These types of algorithms are commonly used in the ETP field \citep{elloumi2014classroom, muller2016real}

Additionally, we introduce a Layer-Cake heuristic algorithm, which builds the exam schedule iteratively using MIP for a subset of exams in a greedy manner. The idea of combining mathematical programming and heuristics appears in many other contexts is what is called ``matheuristics" \citep{archetti2014survey, sorensen2015metaheuristics, boschetti2022matheuristics, siew2024survey}. \cite{gogos2012improved} and \cite{ccimen2022invigilators}  have explored applying matheuristics in general ETPs.

\section{Modeling Approach}\label{sec:model}
This section describes the various components of the MIP-based scheduling frameworks. We begin by introducing the multi-stage Group-then-Sequence model, explaining it stage by stage. Following this, we outline the ``Layer-Cake" heuristic algorithm. Finally, we discuss approaches for creating additional variants of the scheduling models.

\begin{figure}[htbp!]
\centering
\includegraphics[width=\textwidth]{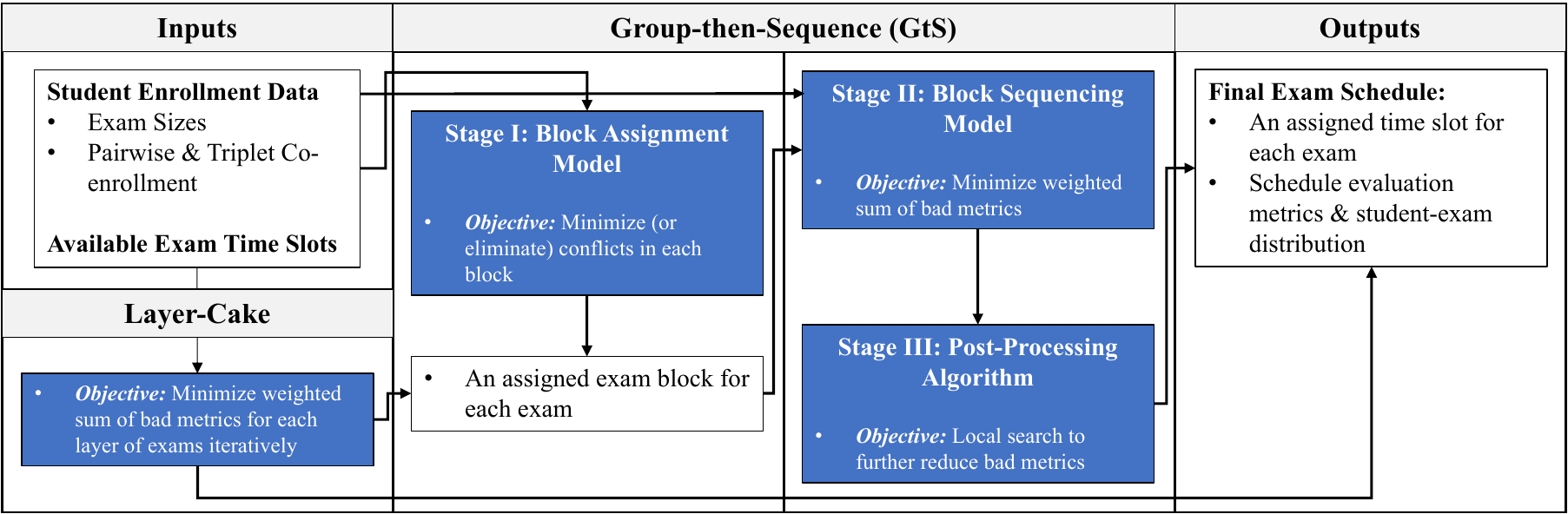}
\caption{Workflow of the Proposed Final Exam Scheduling Framework}
\label{fig:workflow}
\end{figure}

Figure \ref{fig:workflow} illustrates the input and output components of the scheduling framework. The primary source of input data is the student-level enrollment information, provided by the University Registrar in a spreadsheet format each semester. Each row in the spreadsheet represents a student-course pair. To facilitate the scheduling process, this data is further processed to generate three sets of course-level input data required for the subsequent models / algorithms. These include the exam sizes, the number of co-enrolled students for each pair of exams, and the number of co-enrolled students for each triplet of exams. The output of the framework is a mapping that assigns each exam to a specific time slot.

\subsection{Multi-Stage Group-then-Sequence (GtS) Model}
The core of the GtS model consists of the Block Assignment stage and the Block Sequencing stage. A Post-Processing stage can be applied to the schedule generated by GtS. The combined three-stage model is called the GtSP (Group-then-Sequence-Post-Processing).

The idea behind the first two stages is to initially group the exams into a smaller number of exam \textit{blocks} and then assign each block to a time slot. 

\subsubsection{Block Assignment}
The Block Assignment stage is formulated as a MIP problem (\textsc{Min\_Conflict\_IP}) that assigns each exam to one of a specified number of blocks, denoted by $k$. It is assumed that $k \leq |S|$, where $S := \{0, 1, \ldots, |S|-1\}$ is the set of available time slots. The objective of this grouping is to minimize the total number of pairwise co-enrollments (direct conflicts) among exams within the same block. \revision{Appendix \ref{sec:min-ip} provides the \textsc{Min\_Conflict\_IP} formulation.}  \revision{It is important to note that, at this stage, each block represents a group of exams that will be scheduled simultaneously, but is not yet associated with a specific time slot.}

If it is possible to achieve zero conflicts within each block, the block assignment problem can be interpreted as a graph coloring problem. Consider a graph $G = (V, E)$, where $V$ represents the set of exams, and an edge $(i,j) \in E$ connects exams $i$ and $j$ if at least one student is enrolled in both. A proper coloring of the graph $G$ corresponds to a zero-conflict block assignment, as exams with the same color would be grouped in the same block.

We also propose a variation of the Block Assignment problem, named \textsc{Zero\_Conflict\_IP}, that imposes a strict zero-conflict constraint. \revision{We assume a fixed linear ordering of exam blocks $0, 1, \ldots, k-1$ and declare only the consecutive pairs $(b, b+1)$, for $b = 0, \ldots, k-2$, to be neighboring (so block $0$ and block $k-1$ are not neighbors). The model imposes a hard constraint forbidding any co‐enrolled exams from occupying the same block, and it minimizes the total number of co‐enrollments across all such consecutive block pairs. This objective draws an analogy to the Traveling Salesman Problem (TSP)---except that we optimize a Hamiltonian path rather than a cycle, reflecting the linear sequence of blocks. \revision{The \textsc{Zero\_Conflict\_IP} formulation is available in Appendix \ref{sec:zero-ip}.} Figure \ref{fig:block_assignment_viz} shows a three‐block example together with its path‐objective value.} Compared to the minimum-conflict formulation, \textsc{Zero\_Conflict\_IP} has a more constrained feasible space, which increases the likelihood of infeasible solutions. 

The feasibility of having a zero-conflict block assignment depends mainly on the required number of blocks $k$. While one way to find a lower bound on $k$ is to enumerate all possible values and check the feasibility of the \textsc{Zero\_Conflict\_IP} problem, a smarter way is to solve for the Maximum Clique (\textsc{Max\_Clique\_IP}) problem. \revision{Appendix \ref{sec:max_clique_ip} contains the \textsc{Max\_Clique\_IP} formulation.} The \textsc{Max\_Clique\_IP} involves finding a subset of vertices in $G$ that forms a complete sub-graph, where every two distinct vertices are adjacent. In our context, the size of the max clique corresponds to a lower bound on the number of blocks for which it is feasible to have no conflicts.

\begin{figure}[htbp!]
\centering
\includegraphics[width=\textwidth]{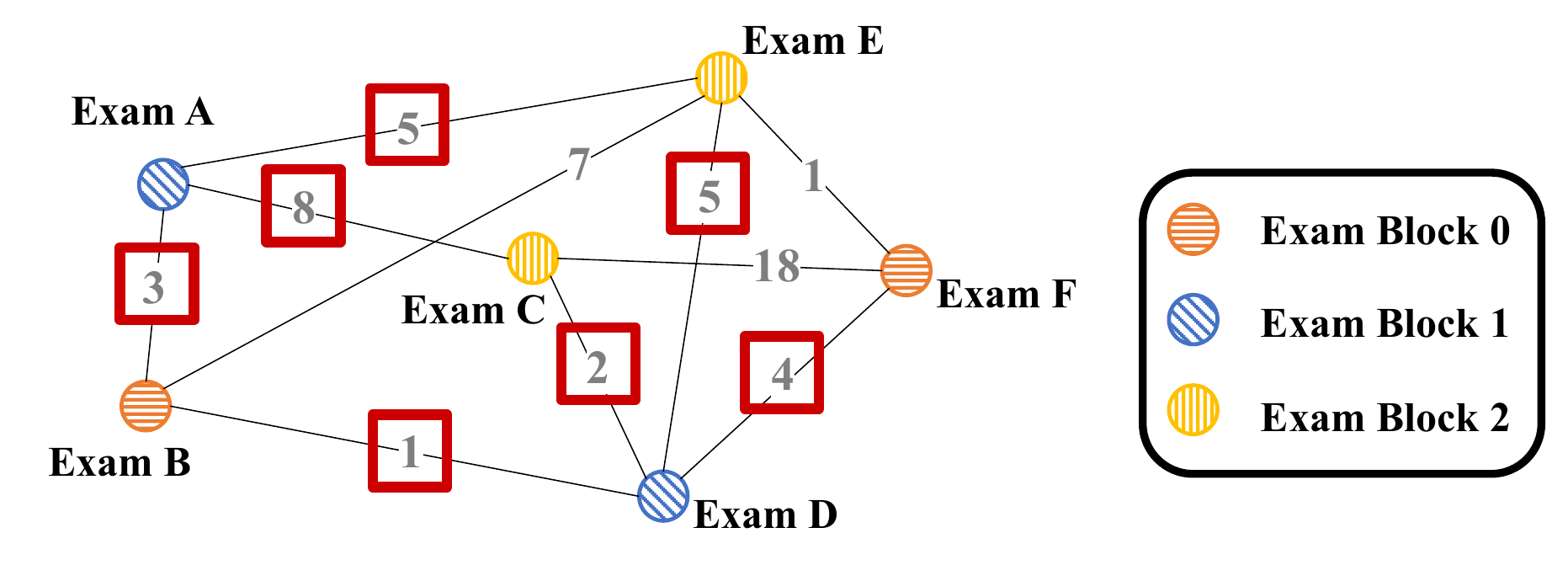}
\caption{This is an example of a zero-conflict block assignment. The value on each edge represents the number of co-enrollments. Blocks 0 and 1, as well as Blocks 1 and 2, are considered neighboring blocks. Consequently, the path-objective value can be calculated as the sum of the highlighted values.}
\label{fig:block_assignment_viz}
\end{figure}

\subsubsection{Block Sequencing}
After assigning exams to the designated blocks, we then allocate these blocks to time slots by \textit{sequencing} them in Block Sequencing stage (\textsc{Block\_Sequencing\_IP}). In the Block Assignment phase, we focus mainly on direct conflicts, but in Block Sequencing, we begin to account for other levels of undesirable metrics. \revision{Appendix \ref{sec:sequencing_ip} presents the \textsc{Block\_Sequencing\_IP} formulation.}

Assume for now that the number of blocks equals the number of time slots, i.e., $|S| = k$. Then, the set of exam blocks is denoted by $B := \{0, 1, \ldots, |S| - 1\}$. To maintain linearity in the objective function, we introduce the four-index variables, denoted as $x_{ij\ell s}$, which allow us to penalize metrics involving triplets of exams. $x_{ij\ell s} = 1$ when block $i$ is assigned to time slot $s$, block $j$ is assigned to time slot $(s + 1) \mod |S|$, and block $\ell$ is assigned to time slot $(s + 2) \mod |S|$. The one-to-one mapping between blocks and time slots necessitates that the four-index variables adhere to a cyclic self-consistent property.

\begin{figure}[htbp!]
\centering
\includegraphics[width=\textwidth]{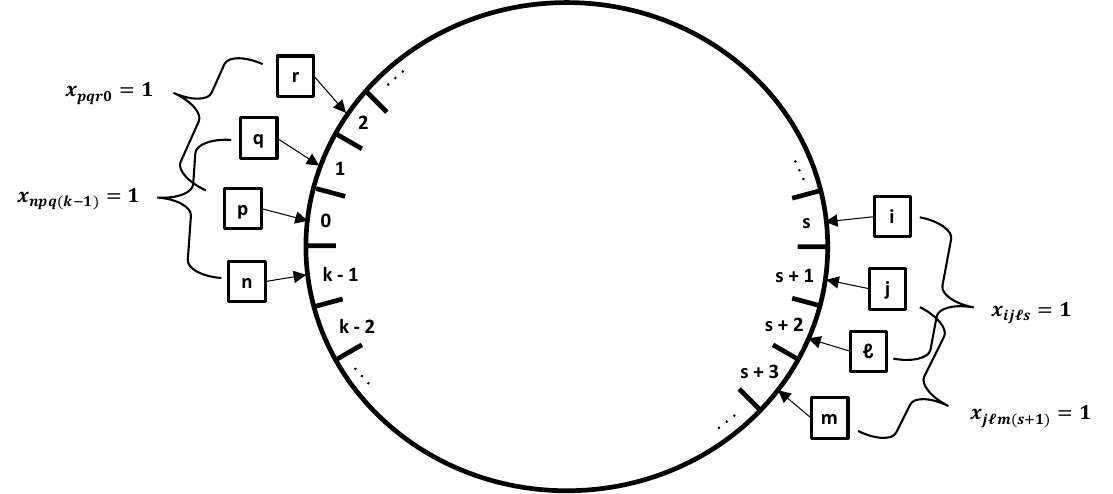}
\caption{Illustration of the $x_{ij\ell s} $ Variables and the Cyclic Self-Consistent Property}
\label{fig:cycle}
\end{figure}

As shown in Figure \ref{fig:cycle}, Block Sequencing produces a mapping of blocks to time slots that are evenly distributed and arranged in a cyclic order. Observe that if $x_{ij\ell s}  = 1$, then there must be another block $m$ assigned to time slot $(s + 3)$, leading to $x_{j\ell m (s+1)} = 1$. This relationship also holds at the end of the time slot sequence due to the use of the modulo function. For instance, if $x_{npq(k-1)} = 1$, then there must be another block $r$ assigned to time slot ($k \mod k = 0$), making $x_{pqr0} = 1$. This observation can be generalized, leading to the formulation of the following self-consistent constraint:
\begin{equation}
    \sum\limits_{i \in B}x_{ij\ell s}  = \sum\limits_{m \in B}x_{j \ell m n_s}, \forall j,\ell \in B, s \in S, \\
\end{equation}
where $n_s:=(s+1) \mod |S|$.

\begin{figure}[htbp!]
\centering
\includegraphics[width=0.4\textwidth]{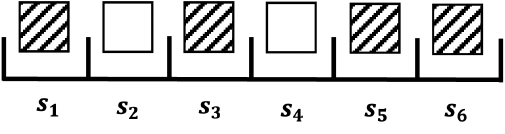}
\caption{An Example Usage of Virtual Blocks}
\label{fig:virtual_block}
\end{figure}

To handle the situation where the number of available time slots exceeds the number of blocks, i.e., $|S| > k$, we introduce an additional layer of optimization. To match the number of blocks with the number of slots, we introduce the concept of a ``virtual block" --- an empty block without any exams. By adding an appropriate number of virtual blocks, we ensure that the total number of blocks equals the number of slots. In the block-to-slot assignment generated by the model, time slots assigned to virtual blocks are considered empty and left unused. This approach helps reduce back-to-back exams and improve other metrics without the need to introduce more blocks. Additionally, virtual blocks enable the exclusion of specific time slots by enforcing constraints that assign virtual blocks to those designated time slots. In Figure \ref{fig:virtual_block}, unshaded squares represent virtual blocks, indicating that slots $s_2$ and $s_4$ are unused in this example.

Additionally, the Block Sequencing model allows for the front-loading of large exams. We introduce hard constraints based on a user-defined exam size cut-off and an early time slot cut-off, ensuring that exams of courses with enrollments exceeding a certain threshold are scheduled earlier in the exam week. Note that in Block Sequencing, the front-loading applies to the entire block containing the large exams rather than scheduling individual exams earlier.

As a result of the Block Sequencing model, each block is allocated to a specific time slot, as are the exams within each block. This produces a working final exam schedule but might not be fully optimized across all metrics. 

\subsubsection{Post-Processing Local Search Algorithm}
Multiple optimal solutions with the same number of direct conflicts can be generated in the Block Assignment stage. However, the schedules derived from these assignments may vary in other metrics, as the Block Assignment model only focuses on direct conflicts. As a result, the block assignment obtained may not be optimal when considering all metrics comprehensively. To address this limitation, we develop a local search algorithm that takes an existing exam schedule and iteratively reschedules subsets of exams with a specified size $n$. \revision{The full pseudocode for the Post-Processing local search algorithm is available in Appendix \ref{sec:post_processing_appendix}.}

In each iteration, the algorithm selects the top $n$ exams to reschedule by assessing the contribution of each exam to the number of bad metrics, normalized by the number of students taking the exam. For example, an exam where 30\% of its students have an exam right before or after it would be prioritized for rescheduling compared to an exam with only 5\% of students in a similar situation. 
\revision{This normalization effectively ranks exams by the density of ``bad" events rather than their absolute counts, allowing the algorithm to focus on the intensity of issues experienced per student. Empirically, this approach provides a more effective signal for guiding the local search.}

Once the subset is selected, the algorithm solves a \textsc{Schedule\_IP} problem to move the exams to time slots that minimize the number of B2Bs and 2-in-24hrs, without increasing the number of direct conflicts. 
The process continues iteratively, \revision{with the algorithm using adaptive step sizes to adjust the starting index of the rolling window with size 
$n$.} If no improvement is found after rescheduling, a larger step forward is taken to explore new solutions. Conversely, if an improvement is observed, a smaller step backward is adopted to fine-tune the schedule. This adaptive approach enables the algorithm to efficiently navigate the solution space, balancing exploration with refinement.

After applying the Post-Processing algorithm to the schedule generated by Block Sequencing, we should obtain a new schedule with significant improvements in reducing the number of bad metrics.

\subsection{Layer-Cake Heuristic Algorithm}

The Layer-Cake heuristic algorithm offers an alternative to the GtS model and is inspired by the Post-Processing algorithm. Unlike the GtS model, Layer-Cake addresses conflicts and other metrics simultaneously by solving the \textsc{Schedule\_IP} problem, which optimizes all metrics together. To reduce the computational complexity of scheduling all exams at once, the problem is solved for a smaller subset of exams each time. A key difference between Layer-Cake and GtS is that Layer-Cake directly takes in a sequence of time slots, eliminating the need to specify the number of blocks, $k$, as an input. \revision{Appendix \ref{sec:layer_cake_appendix} provides the full algorithm pseudocode for Layer-Cake.}

\begin{figure}[ht!]
\centering
\includegraphics[width=0.8\textwidth]{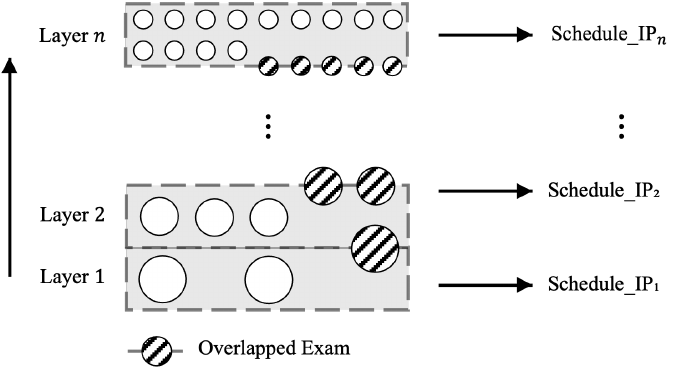}
\caption{\revision{Visualization of the Layer-Cake Algorithm with Partial Overlaps in Consecutive Layers}}
\label{fig:layer_cake}
\end{figure}

Layer-Cake follows a layered approach, where each layer consists of a set of exams such that the total number of students per layer is similar. The layers are determined by sorting the exams in descending order of student enrollment. The initial layers focus primarily on larger exams, while subsequent layers include many smaller exams. The algorithm starts by scheduling the largest exams independently of the others. \revision{It then schedules the next subset of exams to build upon those already assigned, gradually building the schedule ``layer by layer." Crucially, each new layer can partially overlap with the preceding one, as depicted in Figure \ref{fig:layer_cake}.} This ordering ensures that the largest exams are scheduled first, as the \textsc{Schedule\_IP} problem becomes more challenging when dealing with larger exams compared to smaller ones, even with a constant total number of students.

An essential feature of the Layer-Cake algorithm is the overlap between consecutive layers. \revision{In each iteration, a subset of exams that have already been scheduled are included in the IP, allowing them to be rescheduled alongside the newly introduced exams.} This overlapping of layers, serving as a warm start for the IP problems, considerably boosts the algorithm's performance.

Since each \textsc{Schedule\_IP} problem is unaware of exams in future layers and only a subset of the exams can be moved in each iteration, the Layer-Cake algorithm may not yield perfect results. However, it offers notable advantages in computational efficiency over the GtS model. By producing potentially improved schedules in much less time, Layer-Cake provides an effective balance between speed and performance. In addition, the flexibility to adjust the layer size allows for easy customization of this trade-off: larger layers tend to produce better schedules, though at the cost of increased solving time.

\subsection{The Hybrid Model} \label{sec:base-lch}
Building on the strengths of both the GtS model and Layer-Cake, we also develop a new variant of an exam scheduling model. This Hybrid model is inspired by the observation that the assignment of exams to time slots in Layer-Cake, without considering their order, can be interpreted as a block assignment. This opens up the possibility of using the Layer-Cake solution as a promising feasible solution for the Block Assignment model.

Integrating the block assignment obtained from Layer-Cake into the subsequent stages of the GtS model has the potential to improve the overall schedule. Since Layer-Cake considers multiple levels of metrics beyond just direct conflicts, incorporating its block assignment into the GtS model could lead to enhanced schedules compared to the standard Min-Conflict or Zero-Conflict Block Assignment, which primarily focuses on minimizing direct conflicts.

\section{Implementation Results and Impact}\label{sec:implementation}
The proposed final exam scheduling framework has been successfully adopted by the \University\ Registrar for five consecutive semesters, beginning in Spring 2022. The scheduling process usually starts about two weeks after the semester begins, allowing the Registrar to collect final enrollment information during the Add/Drop period. To efficiently create the final exam schedule within a tight timeline of less than a month, a dedicated team of undergraduate students is trained in advance.

Once the team gains access to the necessary data, we conduct experiments using various parameters and models implemented in Gurobi and Python. This collaborative approach allows for comprehensive exploration and evaluation of different scheduling options. To accelerate the process, we use Google Colab to run multiple experiments in parallel. From the generated subset of ``good" schedules, the University Registrar selects the best one to be published on the university website.

In the subsequent parts, we showcase the performance of our scheduling models on real student enrollment data, highlighting the framework's ability to outperform the historical baseline and its flexibility in meeting the registrar's needs. \revision{All experiments in this paper are conducted on a single machine equipped with an 8-core Intel Xeon CPU running at 2.20 GHz. All optimization models are solved using Gurobi version 11.0.3. Computation times are not reported in this paper, as they are not a primary concern in our application. All methods complete within several hours, which is acceptable given that the final exam schedule is typically prepared over the span of approximately one month each semester.}

\subsection{Performance of the Scheduling Models Across Multiple Semesters}

In the first set of experiments, we compare the performance of different scheduling models using actual student enrollment data at \University\ from Spring 2022 to Spring 2024. To support open-source research, an anonymized version of Cornell’s Spring 2024 exam-scheduling input data has been made publicly available at \href{https://github.com/AdamJovine/SchedulingArchiveCornell}{https://github.com/AdamJovine/SchedulingArchiveCornell}.
Table \ref{tab:stats} presents descriptive statistics for student enrollment and exam timing data. Generally, more students take exams in the fall semester than in the spring. The last column shows the Max-Clique size, representing a lower bound on the number of exam blocks required for that semester. For instance, in Spring 2024, a Max-Clique of size 20 indicates that there was a subset of 20 courses with overlapping enrollments, necessitating at least 20 distinct exam blocks to avoid any direct conflicts.

\begin{table}[htbp!]\scriptsize
\centering
\caption{Descriptive Statistics of Input Data From Spring 2022 to Spring 2024 }
\label{tab:stats}
\begin{tabular}{ c c c c c c }
\hline
 & \# of Students & \revision{\# of Unique Exams} & Exam Time Range & \# of Available Time Slots & Max-Clique Size\\
\hline
SP24 & 16265 & 544 & May 11th 9 a.m. - May 18th 7 p.m. & 24 (24) & 20 \\
FA23 & 18319 & 601 & Dec. 8th 9 a.m. - Dec. 15th 7 .m & 24 (23) & 18 \\
SP22 & 16261 & 553 & May 14th 9 a.m. - May 21th 7 p.m. & 24 (19) & 17 \\
FA22 & 18081 & 588 & Dec. 9th 9 a.m. - Dec. 17th 9 a.m. & 25 (22) & 19\\
SP22 &  16429 & 539 & May 13th 9 a.m. - May 20th 7 p.m. & 24 (20)$^{*}$  & 16\\ 
\hline
\end{tabular}
\vspace{6pt} 
    
    \raggedright
    Note: *The numbers inside the parentheses correspond to the actual number of available time slots, which is constrained by the University Registrar but subject to relaxation starting from SP24.
\end{table}

\begin{table}[htbp!]
    \scriptsize
    \centering
    \caption{Performance of the Baseline and Our Scheduling Models Across 5 Semesters}
    \label{tab:compare-semester-part1}
\begin{tabular}{l|lccccc}
        \toprule
        Semester & Model & Conflict & Triple & B2B & 2-in-24hr & 3-in-4 \\
        \midrule
        SP24 & Baseline & 790 & 128 & 3029 & 2653 & 614 \\
             & GtS (Ours) & 0 & 50 & 2324 & 2790 & 441 \\
             & GtSP (Ours) & 0 & 30 & 1753 & 2474 & 285 \\
             & Zero-GtSP (Ours) & 0 & 39 & 1763 & 2884 & 357 \\
             & Layer-Cake (Ours) & 13 & 26 & 2020 & 2531 & 406 \\
             & Hybrid (Ours) & 6 & 25 & 1812 & 2236 & 325 \\

        \midrule
        FA23 & Baseline & 376 & 143 & 2883 & 3930 & 736 \\
             & GtS (Ours) & 0 & 51 & 2938 & 3534 & 641 \\
             & GtSP (Ours) & 0 & 48 & 2215 & 3669 & 490 \\
             & Zero-GtSP (Ours) & 0 & 33 & 2038 & 3185 & 551 \\
             & Layer-Cake (Ours) & 0 & 26 & 2323 & 3088 & 546 \\
             & Hybrid (Ours) & 0 & 26 & 1950 & 2775 & 458 \\ 

        \midrule
        SP23 & Baseline & 720 & 188 & 2957 & 2600 & 742 \\
             & GtS (Ours) & 0 & 43 & 2332 & 2709 & 430 \\
             & GtSP (Ours) & 0 & 40 & 2024 & 2211 & 281 \\
             & Zero-GtSP (Ours) & 0 & 27 & 1651 & 2648 & 339 \\
             & Layer-Cake (Ours) & 0 & 22 & 1784 & 2545 & 420 \\
             & Hybrid (Ours) & 0 & 22 & 1566 & 2309 & 409 \\

        \midrule
        FA22 & Baseline & 614 & 142 & 3284 & 3930 & 812 \\
             & GtS (Ours) & 0 & 31 & 2210 & 3556 & 667 \\
             & GtSP (Ours) & 0 & 32 & 1888 & 2942 & 541 \\
             & Zero-GtSP (Ours) & 0 & 37 & 2051 & 2625 & 466 \\
             & Layer-Cake (Ours) & 1 & 25 & 2233 & 3100 & 427 \\
             & Hybrid (Ours) & 1 & 34 & 1621 & 2726 & 299 \\

        \midrule
        SP22 & Baseline & 564 & 206 & 3050 & 2801 & 564 \\
             & GtS (Ours) & 0 & 53 & 2494 & 2493 & 380 \\
             & GtSP (Ours) & 0 & 28 & 1753 & 2309 & 285 \\
             & Zero-GtSP (Ours) & 0 & 44 & 1712 & 2941 & 466 \\
             & Layer-Cake (Ours) & 0 & 20 & 1610 & 2782 & 341 \\
             & Hybrid (Ours) & 0 & 16 & 1411 & 2517 & 264 \\
        \midrule
        Average & Baseline & 612.8 & 161.4 & 3040.6 & 3182.8 & 693.6 \\
        & GtS (Ours) & 0 & 45.6 & 2459.6 & 3016.4 & 511.8 \\
        & GtSP (Ours) & 0 & 35.6 & 1926.6 & 2721 & 376.4 \\
        & Zero-GtSP (Ours) & 0 & 36 & 1843 & 2856.6 & 435.8 \\
        & Layer-Cake (Ours) & 2.8 & 23.8 & 1994 & 2809.2 & 428 \\
        & Hybrid (Ours) & 1.4 & 24.6 & 1672 & 2512.6 & 351 \\

        \bottomrule
    \end{tabular}

\end{table}

To ensure a fair comparison, we use the full set of time slots for all models as shown in Table \ref{tab:stats}. \revision{For each semester, the number of blocks used remains consistent across different scheduling model variants.} For Layer-Cake-based models, the time slot configuration input matches that obtained from the corresponding GtS models. Unless otherwise specified, we use the standard scheduling parameters given in Appendix \ref{sec:setting_appendix} throughout the paper. The Baseline model follows the historical approach described in Section \ref{sec:intro}. The current course schedule includes new meeting time patterns that are not accounted for in the lecture time-based mapping used in the baseline approach. To adapt the Baseline model for different semesters, we assigned courses with new meeting times greedily to one of the existing time slots. This approach preserved the direct mapping between meeting times and exam times. 

Table \ref{tab:compare-semester-part1} summarizes the metrics of the schedules created by our framework and the Baseline. The following key observations can be drawn from the table.

\subsubsection{Dominance Over the Historical Approach}

The Baseline model, commonly used in many universities, is completely dominated by our scheduling models every semester. 
The limitations of relying solely on lecture times to create final exam schedules are evident, as the Baseline schedules consistently underperform across nearly all metrics, particularly in reducing direct conflicts and triples. On average, the Baseline results in 612.8 direct conflicts and 161.4 triples, leading to over 774 occurrences of exam rescheduling. In contrast, schedules generated by our framework generally have zero or minimal conflicts and fewer than 50 triples.

\subsubsection{Impact of Post-Processing}

Comparing the GtS and GtSP models reveals that the Post-Processing stage significantly improves all metrics for the GtS model. This underscores the limitations of the draft schedule generated in the initial two stages of the GtS model, which focus on block-level metrics rather than individual exams. On average, we observe reductions of 22\% in B2Bs, 22\% in triples, 10\% in 2-in-24hrs, and 26\% in 3-in-4s from GtS to GtSP .

\subsubsection{GtSP vs. Layer-Cake}
The Layer-Cake heuristic algorithm presents a distinct trade-off among metrics compared to the GtS-based models, primarily due to its slightly different objective.
Layer-Cake generally excels in providing finer control over bad metrics as its \textsc{Schedule\_IP} directly penalizes the occurrence of B2Bs and 2-in-24hrs. As a result, triples and 3-in-4s are indirectly addressed within each \textsc{Schedule\_IP}, since both of these metrics include at least one 2-in-24hr. Table \ref{tab:compare-semester-part1} reveals that, on average, Layer-Cake produces 33\% fewer triples than GtSP. However, this comes at the cost of 3\% more B2Bs and 2-in-24hrs, as well as a 13\% increase in 3-in-4s.

\subsubsection{Layer-Cake vs. Hybrid}

We find that the Hybrid model outperforms Layer-Cake in nearly all metrics, with the exception of triples. On average, compared to Layer-Cake, the Hybrid model reduces B2Bs by 16\%, 2-in-24hrs by 10\%, and 3-in-4s by 18\%. Additionally, the Hybrid model achieves the same or fewer direct conflicts than Layer-Cake, as it is built upon the block assignments generated by the Layer-Cake schedule. The Hybrid model also benefits from the Post-Processing stage, which can further reduce direct conflicts through rescheduling, as demonstrated in SP24. For these reasons, we typically use Layer-Cake for block assignment but rely on the Hybrid model for the final schedule.

\subsubsection{Impact of Different Block Assignment Methods}

Table \ref{tab:compare-semester-part1} also compares different approaches to Block Assignment: \textsc{Min\_Conflict\_IP} (the default GtSP), \textsc{Zero\_Conflict\_IP} (the Zero-GtSP), and using the Layer-Cake schedule for block assignment (the Hybrid).

Among the two IP methods, solving the \textsc{Min\_Conflict\_IP} to optimality is relatively easier. With $k = 22$, the \textsc{Min\_Conflict\_IP} yields a zero-conflict block assignment within an hour.  However, it's important to note that the \textsc{Min\_Conflict\_IP} does not always guarantee zero conflict, particularly with smaller values of $k$. In contrast, the \textsc{Zero\_Conflict\_IP} model guarantees zero conflicts but is more challenging to solve, making it difficult to obtain an optimal or near-optimal solution within a reasonable time limit. When a solution is feasible, the \textsc{Zero\_Conflict\_IP} generally outperforms the \textsc{Min\_Conflict\_IP} in terms of B2Bs due to its TSP objective, as discussed earlier. Specifically, the \textsc{Zero\_Conflict\_IP} reduces B2Bs by 4\% compared to the \textsc{Min\_Conflict\_IP}, but this comes with a 5\% increase in 2-in-24hrs and a 15\% increase in 3-in-4s.

Compared to the IP methods, the Hybrid model leverages the advantages of Layer-Cake by considering higher-level metrics during block assignment. In some semesters, the Hybrid model even outperforms GtSP and Zero-GtSP, while offering a strong alternative in others. On average, the Hybrid model results in 15\% fewer B2Bs and 20\% fewer triples compared to GtSP.

\subsection{Flexibility of the Framework}

The results presented in the following subsections highlight the flexibility of our framework. These results are derived from routine experiments conducted each semester, allowing us to explore a variety of \revision{high-quality schedules} and present trade-offs to the University Registrar. Specifically, we focus on Spring 2024, the most recent semester at the time this paper is written.

\subsubsection{Flexibility in Controlling the Number of Blocks}

Understanding the relationship between the number of blocks and performance of the generated schedule is complex. On one hand, a higher value of $k$  provides more time slots for exam assignments, allowing them to be spread out and potentially improving the metrics. On the other hand, the ``virtual block" design enables strategic placement of empty time slots between blocks with high co-enrollment. 

As shown in Table \ref{tab:experiment-num-block}, the B2B, 2-in-24hr, and 3-in-4 metrics do not follow a linear trend as $k$ increases.  This experiment demonstrates that high-quality schedules can be produced with relatively few blocks, though this was achieved with complete flexibility in the placement of empty time slots.
\begin{table*}[ht!] \scriptsize
    \centering
    \caption{Performance of GtSP with Different Numbers of Blocks}
    \label{tab:experiment-num-block}
    \begin{tabular}{lcccccc} 
        \toprule
        \# Blocks ($k$) & Conflict & Triple & B2B & 2-in-24hr  & 3-in-4\\
        \midrule
        20 & 3 & 17 & 2011 & 2397  & 396 \\
        21 & 0 & 22  & 1950 & 2574 & 286 \\
        22 & 0 & 43 & 1838 & 2487  & 358 \\
        23 & 0 & 30 & 2009 & 2684  & 397 \\
        24 & 0 & 60 & 1910 & 2477  & 283 \\
        
        \bottomrule
    \end{tabular}
\end{table*}

\subsubsection{Flexibility in Restricting Time Slot Availability}

As previously noted, the University Registrar often requests the exclusion of specific time slots to help reduce student stress and accommodate graduation events. To assess the potential impact of these exclusions, we conduct experiments by rerunning Block Sequencing on a 20-block assignment, excluding between 0 and 4 time slots using the virtual block design in Block Sequencing. The  time slots considered for exclusion include a Saturday evening (slot 3), a Sunday morning (slot 4), a Friday evening (slot 21), and a Sunday evening (slot 24).

Table \ref{tab:experiment-excluded-times} presents the results of these exclusions. As more slots are restricted, the model's options to optimize metrics diminishes. We observed a noticeable increase in the number of 2-in-24hrs and triples, suggesting that restricting multiple time slots can result in \revision{less favorable student schedules}. In response, the team has been working with the Registrar to open additional time slots, which take effects starting from SP24.

\begin{table*}[!ht] \scriptsize
    \centering
    \caption{Performance of GtSP with Different Excluded Time Slots}
    \label{tab:experiment-excluded-times}
    \begin{tabular}{lcccccc}
        \toprule
        Excluded Times & Conflict & Triple & B2B & 2-in-24hr  & 3-in-4\\
        \midrule
        
        \{\} & 3 & 17 & 2011 & 2397  & 396 \\
        \{3\}  & 3 & 25 & 2101 & 2297  & 416\\
        \{3, 4\} & 3 & 30 & 1881 & 2573  & 424\\
        \{3, 4, 21\} & 3 & 33 & 1868 & 2474  & 343\\
        \{3, 4, 21, 24\} & 3 & 41 & 1993 & 2772 & 397 \\
        \bottomrule
    \end{tabular}
\end{table*}

\subsubsection{Flexibility in Front-Loading Large Exams} 
In this section, we explore the impact of different exam size and time slot cutoffs for front-loading large exams. \revision{We modify the Block Sequencing stage of the GtSP model by experimenting with a range of size thresholds for large exams and the corresponding early time slot cutoffs.}

\begin{figure}[htbp!]
    \centering
    \begin{subfigure}[b]{0.45\textwidth}
        \centering
        \includegraphics[width=\textwidth]{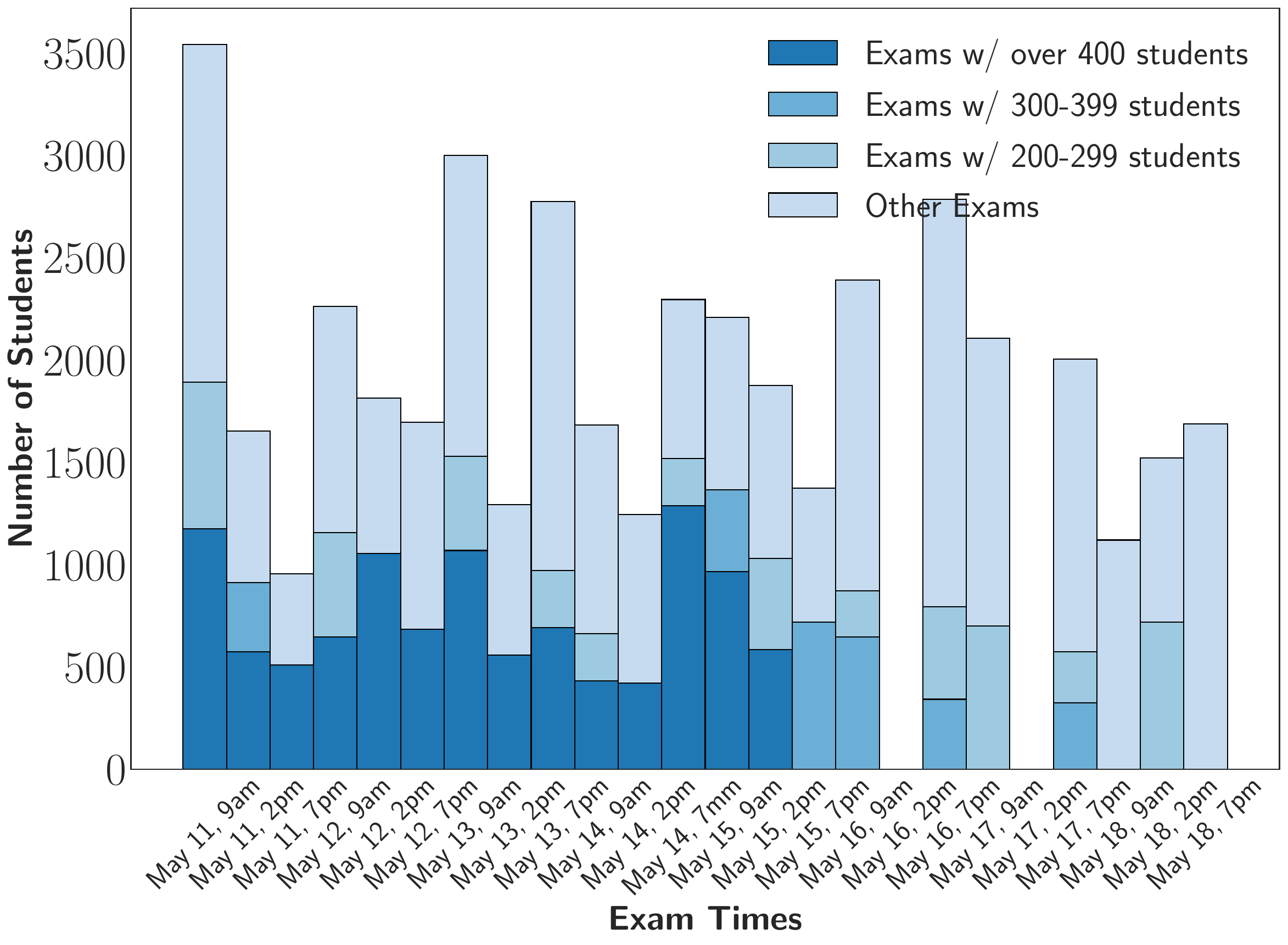}
        \caption{Size Cutoff: 400; Slot Cutoff: 15}
    \end{subfigure}
    \hspace{1em}
    \begin{subfigure}[b]{0.45\textwidth}  
        \centering
        \includegraphics[width=\textwidth]{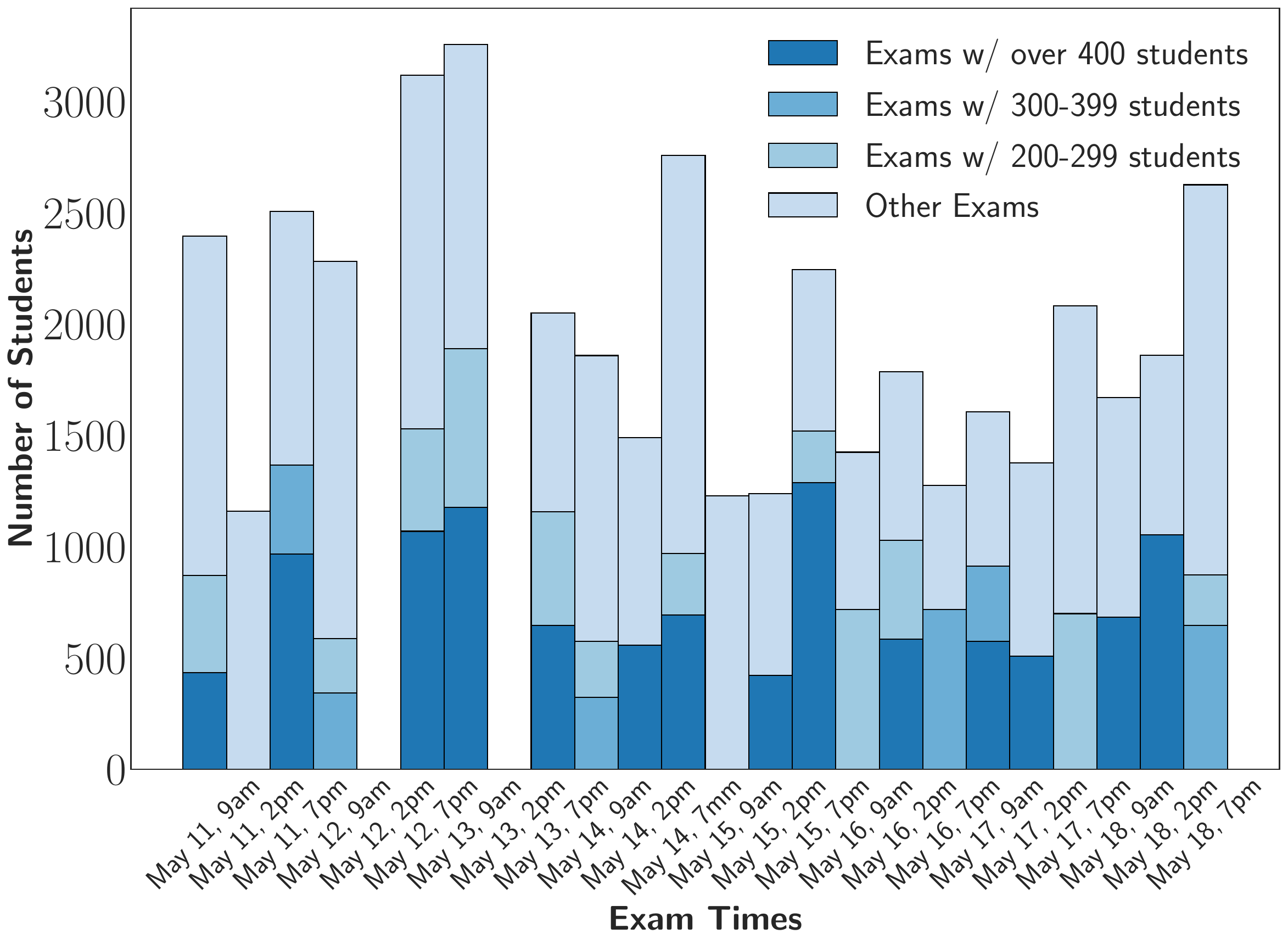}
        \caption{Size Cutoff: 400; Slot Cutoff: 24}
    \end{subfigure}
    
    \vspace{1em} 
    
    \begin{subfigure}[b]{0.45\textwidth}
        \centering
        \includegraphics[width=\textwidth]{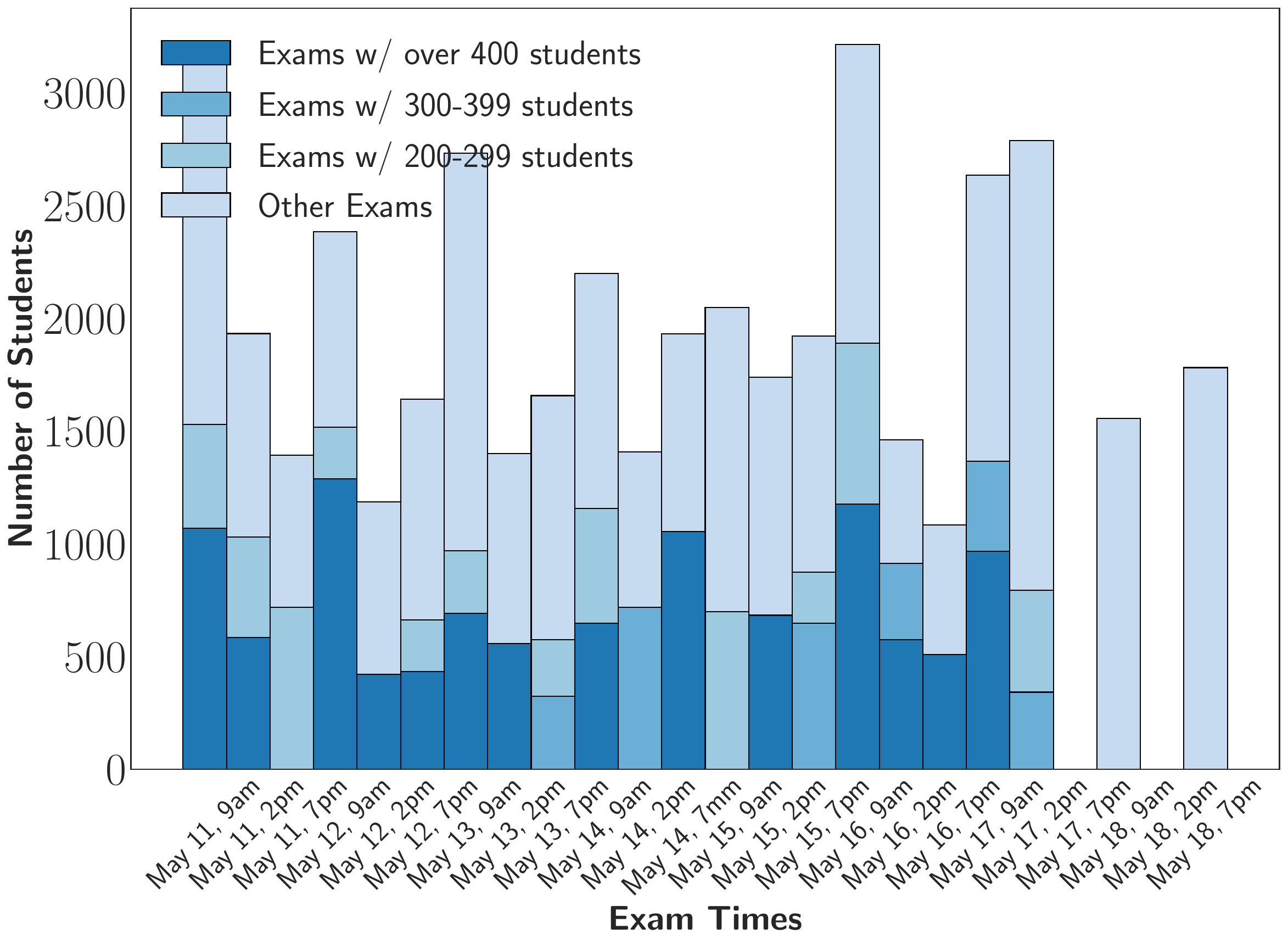}
        \caption{Size Cutoff: 200; Slot Cutoff: 21}
    \end{subfigure}
    \hspace{1em}
    \begin{subfigure}[b]{0.45\textwidth}  
        \centering
        \includegraphics[width=\textwidth]{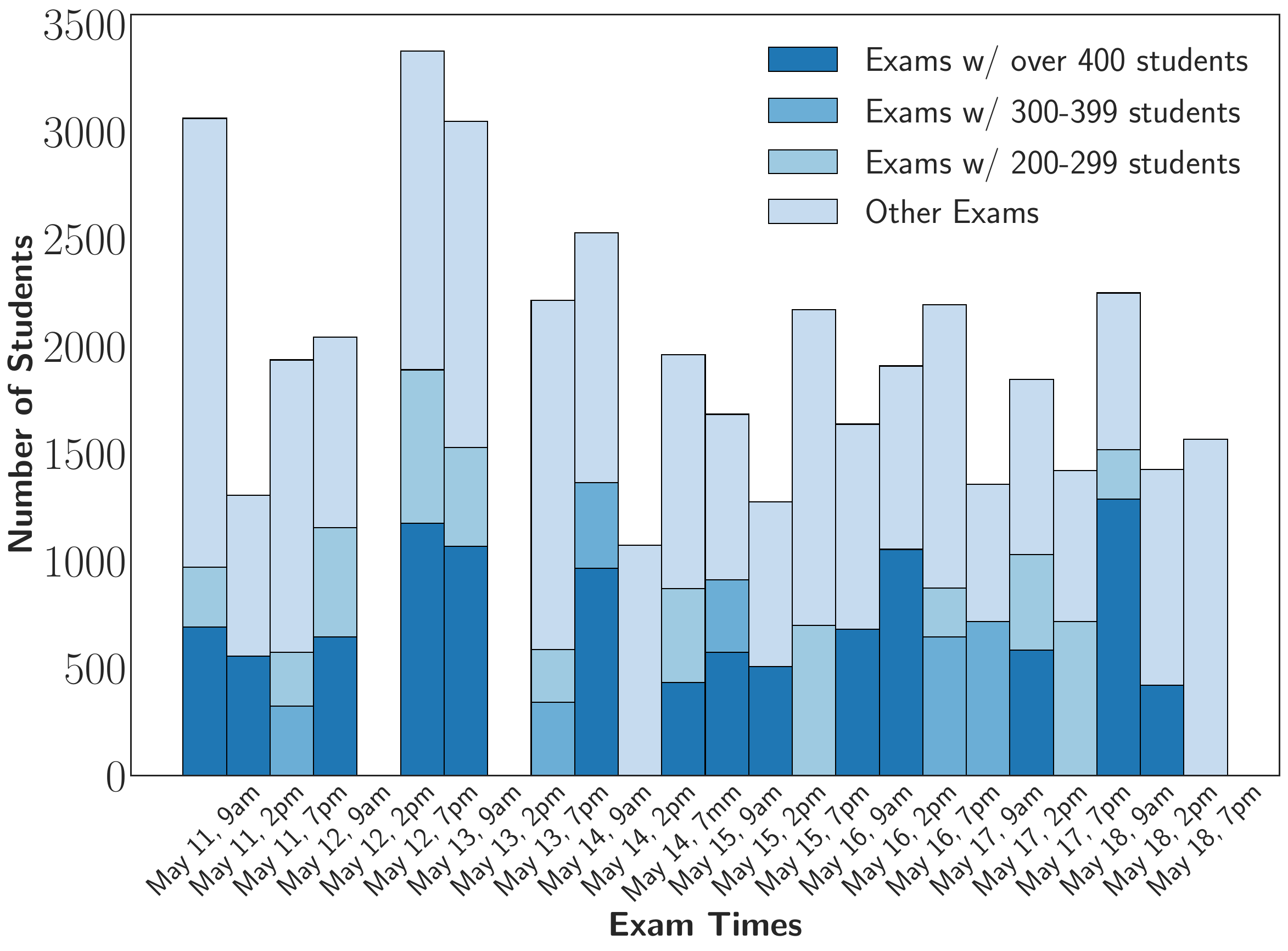}
        \caption{Size Cutoff: 200; Slot Cutoff: 24}
    \end{subfigure}
    \caption{Distribution of Exams with Different Size and Slot Cutoffs}
    \label{fig:frontload-threshold}
\end{figure}

The effects of varying front-loading thresholds align with our expectations, as shown in Figure \ref{fig:frontload-threshold}. Intuitively, the smaller the size cutoff and the earlier the slot cutoff, the more pronounced the front-loading effect becomes. As the size cutoff decreases from 400 to 200 and the slot cutoff moves forward from slot 24, large exams becomes more concentrated in earlier time slots, with only a few scattered large exams remaining on the last few days.  The four sub-figures provided below are visualizations we would include when sharing the final schedules with the University Registrar. These visuals clearly depict the temporal distribution of exams, helping the registrar's office to understand room demand and allocate exams accordingly.

\begin{figure}[htbp!]
    \centering
    \includegraphics[width=0.5\textwidth]{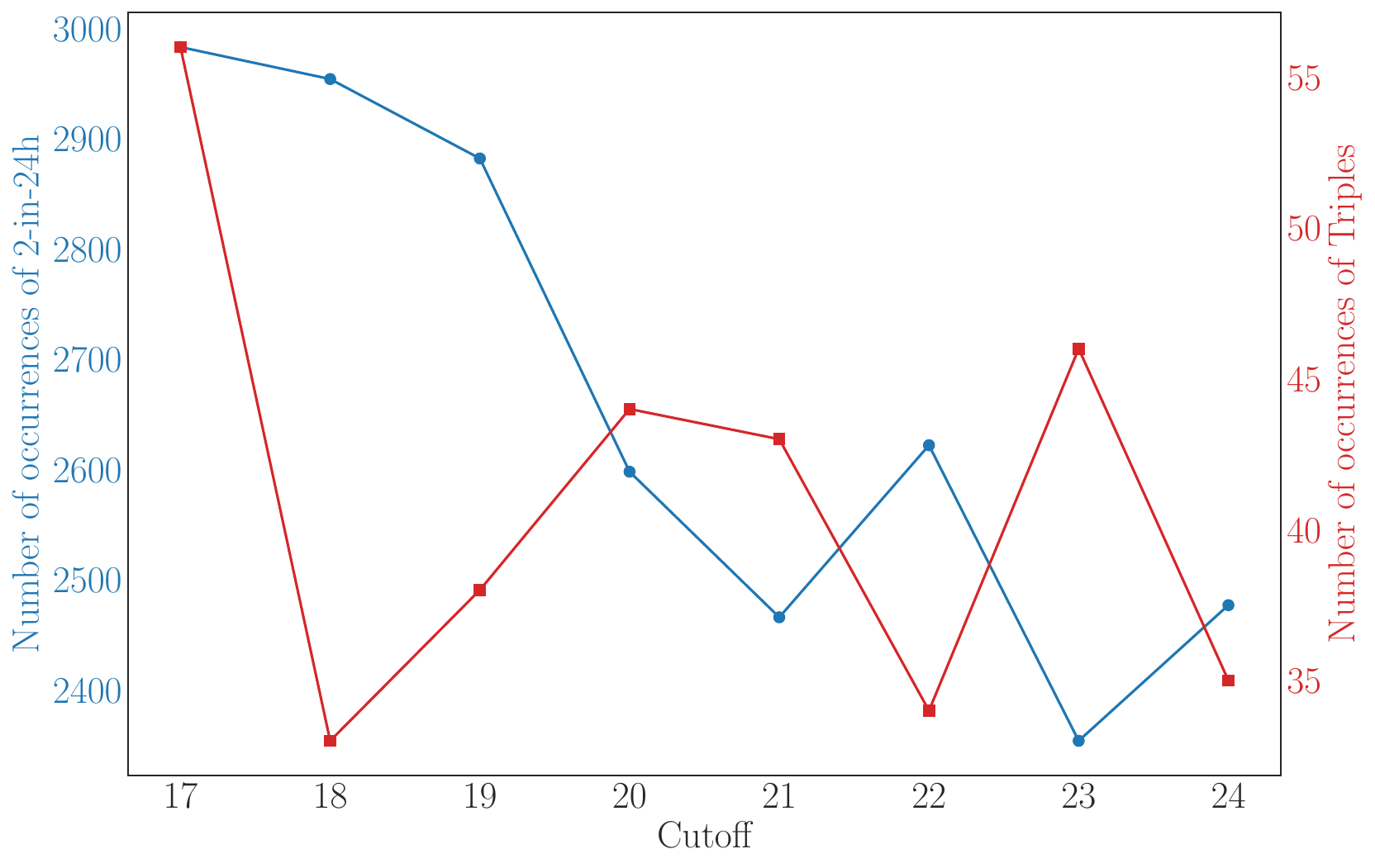}
    \caption{Impact of Different Slot Cutoffs on 2-in-24hrs and Triples}
    \label{fig:slot-cutoff}
\end{figure}

The impact of changing the time slot cutoff on scheduling metrics is further examined and displayed in Figure \ref{fig:slot-cutoff}.
We observe that setting an earlier slot cutoff for large exams is associated with an increase in the number of 2-in-24hrs and triples. This outcome is anticipated, as limiting the available slots for blocks containing large exams reduces the feasible solution space for both Block Sequencing and Post-Processing.

\subsection{Regression Analysis on Key Input Parameters}

In the previous section, we have observed that the number of blocks, metric weights, size cutoff, and slot cutoff can crucially impact the metrics.
\revision{To better understand the trade-offs among these parameters, we use the GtSP model to generate schedules for all combinations of parameter settings listed in Table~\ref{tab:parameters}, evaluated over the three spring semesters.} In total, we produce 2430 distinct schedules (810 per semester) and computed their metrics. We then fit a Random Forest Regression model on these results and perform a SHAP Value analysis \citep{NIPS2017_7062} to evaluate the influence of each parameter on the metrics. SHAP values provide a measure of feature importance, quantifying the contribution of each parameter to the model's output. In this case, we choose the model output to be the number of students who would need rescheduled exams, as this is the most key concern. Recall that, at \University, students with more than two exams in a 24-hour period are eligible for a rescheduled exam.

\begin{table}[!ht]\footnotesize
\centering
\caption{Parameters and Their Values Used for the Regression Analysis}
\label{tab:parameters}
\begin{tabular}{lc}
\toprule
Parameter           & Values                         \\ \midrule
Number of Blocks    & [19, 20, 21, 22, 23, 24]       \\
Size Cutoff         & [150, 200, 250, 300, 350]      \\
Slot Cutoff*        & [16, \ldots, 24] \\
Triple vs. B2B         & [1, 2, 10]                     \\ \bottomrule
\end{tabular}
\vspace{6pt} 
    
\raggedright
\textit{*Note: Schedules were excluded from the data if using a slot cutoff proved infeasible.}
\end{table}

\begin{figure}[htbp!]
    \centering
    \begin{subfigure}[b]{\textwidth}
        \centering  
        \includegraphics[width=0.8\textwidth]{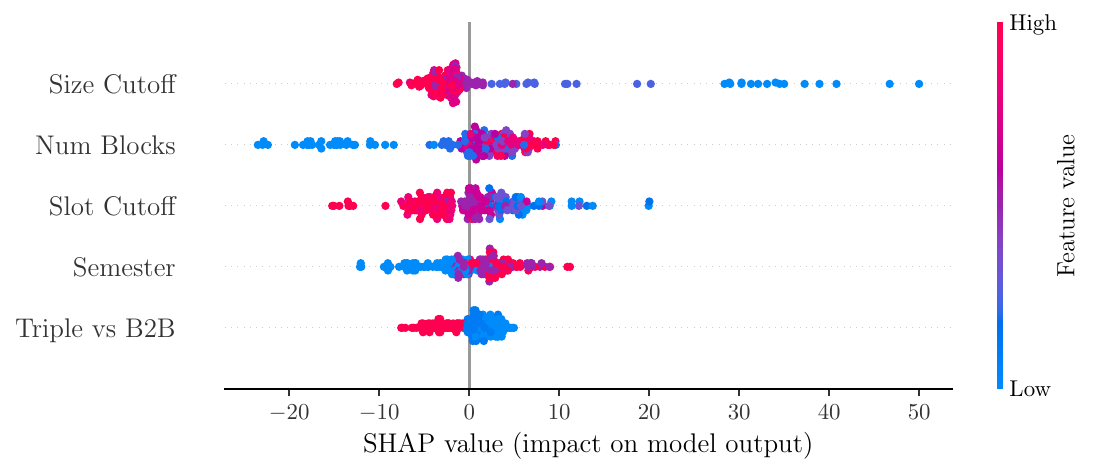}
        \caption{Feature Importance Analysis Using SHAP Values}
    \end{subfigure}
    
    \vspace{1em} 
    
    \begin{subfigure}[b]{0.45\textwidth}
        \centering
        \includegraphics[width=\textwidth]{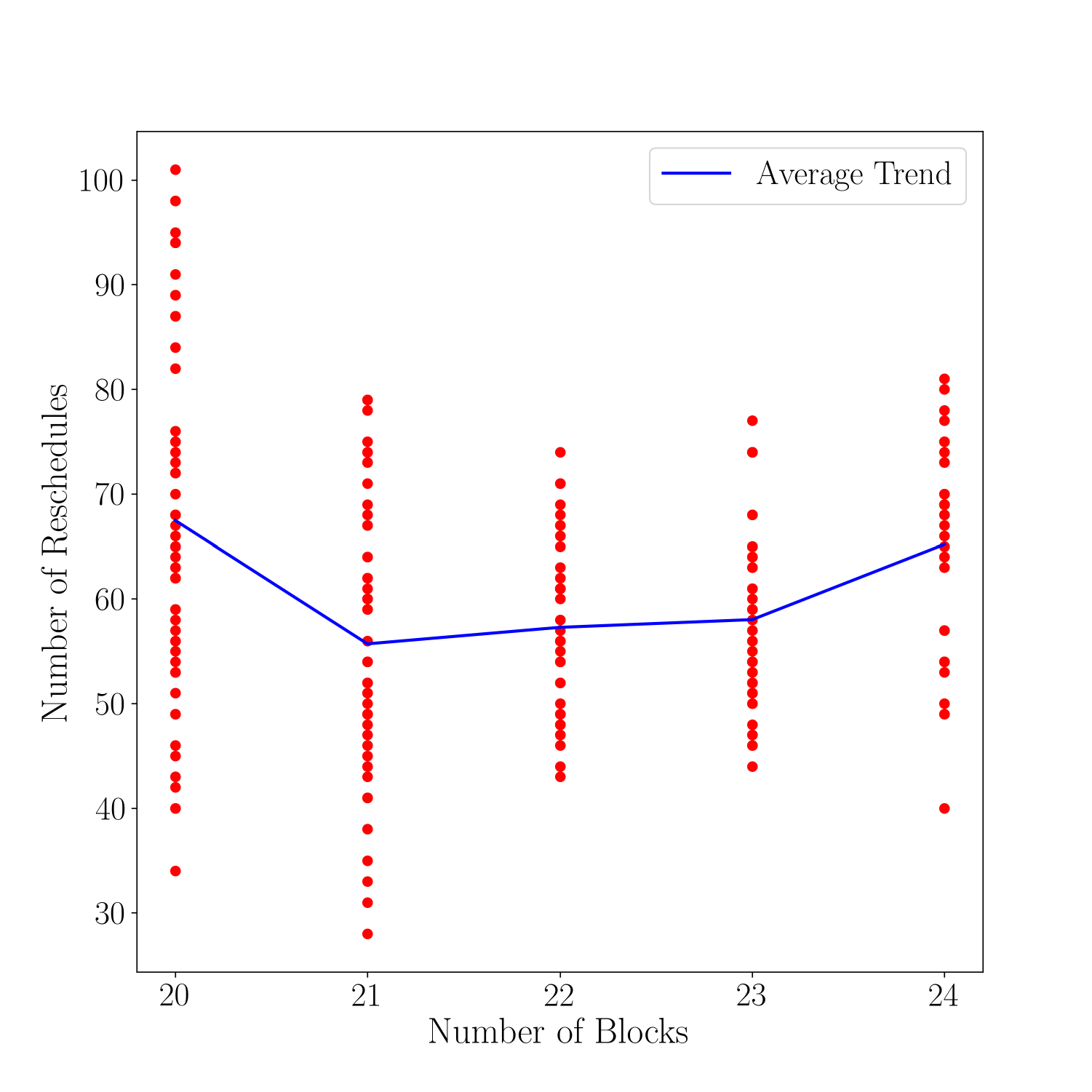}
        \caption{Number of Blocks Vs. Reschedules}
    \end{subfigure}
    \hspace{1em} 
    \begin{subfigure}[b]{0.45\textwidth}
        \centering
        \includegraphics[width=\textwidth]{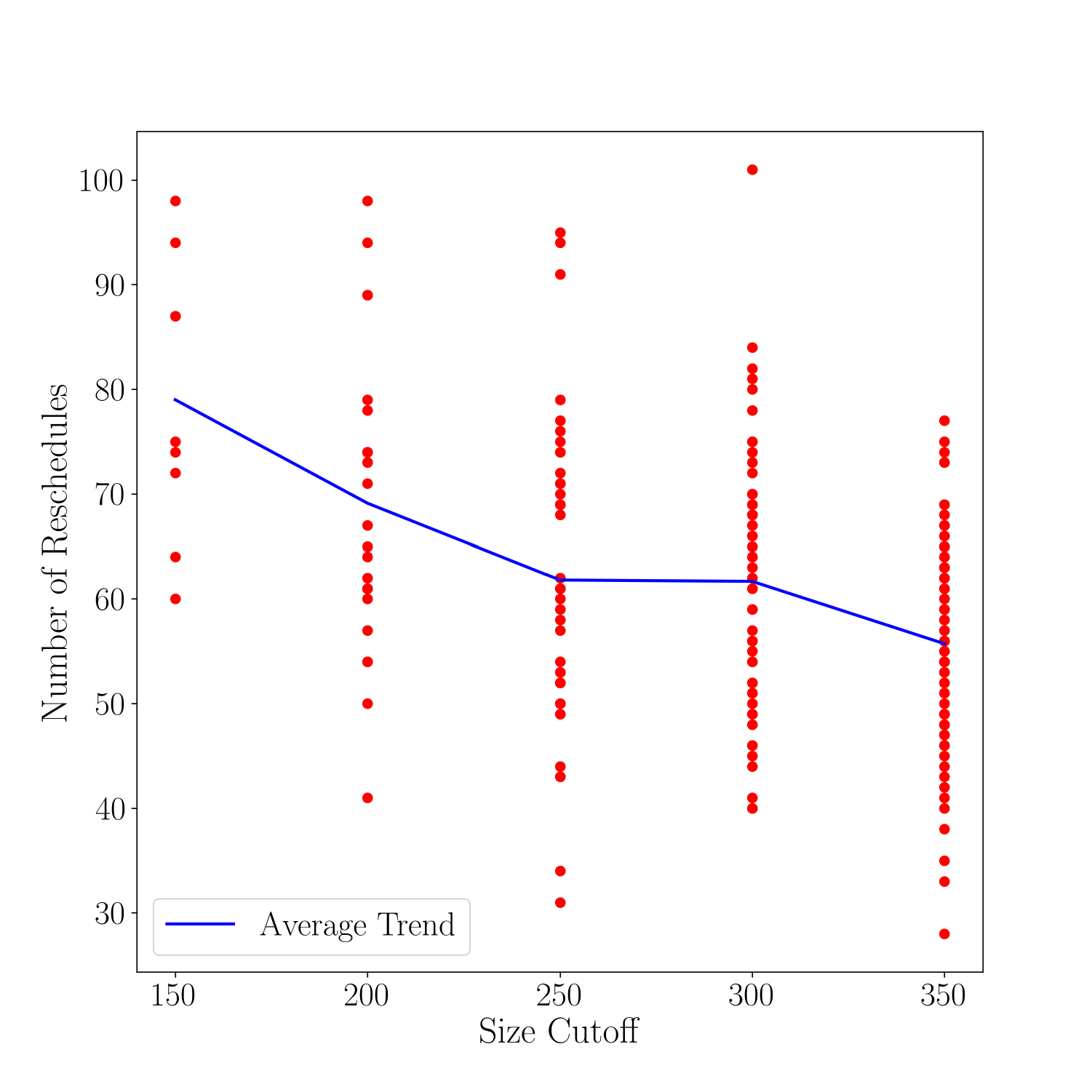}
        \caption{Size Cutoff Vs. Reschedules}
    \end{subfigure}
    
    \caption{Impact of Different Model Parameters on the Number of Rescheduled Exams}
    \label{fig:reschedule-analysis}
\end{figure}
Figure \ref{fig:reschedule-analysis} explores how size cutoffs, time slot cutoffs, and the number of blocks influence scheduling flexibility and the need for rescheduling. Higher size cutoffs generally correlate with fewer reschedules, while low size cutoffs often require higher slot cutoffs to maintain flexibility.

Increasing the number of blocks typically improves schedule performance, though outcomes vary more widely when block numbers are low. Additionally, a higher triple vs B2B trade-off appears to slightly reduce rescheduled exams.  Lastly, schedules from later semesters tend to experience more rescheduling, possibly due to higher enrollment levels.

\section{Conclusions}\label{sec:conclsion}
This paper presents a MIP-based optimization framework that leverages a multi-stage Group-then-Sequence approach combined with a Layer-Cake heuristic algorithm to automate the final exam scheduling process at \University. Over five consecutive semesters, we have collaborated closely with the University Registrar to schedule final exams for the entire student body, consistently generating optimized schedules that outperform historical baselines. \revision{Our results demonstrate not only the effectiveness of our method but also its capacity to adapt to a wide range of constraints and objectives.} \revision{Further evidence of the framework’s generalizability is provided in Appendix~\ref{sec:Nottingham_appendix}, which demonstrates its competitive performance against established exam timetabling algorithms on the University of Nottingham dataset.}

Furthermore, our ongoing partnerships with multiple universities underscore the inherent flexibility of our framework. By accommodating diverse institutional preferences, our approach facilitates the creation of customized exam schedules that balance improved student outcomes with other operational priorities. For instance, if an institution prefers to avoid Sunday morning exams, our model can be readily adjusted to enforce that constraint. Similarly, if minimizing instances of triple exams is a priority, the objective function can be reweighted accordingly. These features highlight the robustness and adaptability of our framework, making it a valuable tool for addressing the complex and varied demands of exam scheduling in higher education.

\ACKNOWLEDGMENT{%
We would like to extend our gratitude to the Cornell ORIE scheduling team members, past and present, who contributed to the development of new models and implementation of our experiments. 

We also express our appreciation to Casey Washburn for providing us with data access and engaging in fruitful discussions. Furthermore, we would like to acknowledge the contributions of Shijin Rajakrishnan, Bob Bland, Dmitriy Levchenkov, and Gwen Spencer in the early stages of the project. }

\bibliographystyle{informs2014} 
\bibliography{refs} 

\clearpage

\begin{APPENDICES}

\section{The Block Assignment Problem} \label{sec:block_assignment_appendix}

\begin{table*}[!hp]\scriptsize
    \centering
    \caption{Notations for the Block Assignment Problem}
    \label{tab:param1}
    \begin{tabular}{p{3cm}p{12cm}p{1cm}p{6cm}}
        \toprule
        \textbf{Notation} & \textbf{Description} \\
        \midrule
        \underline{Sets}  \\
        {$B$} & Set of exam blocks = $\{0, \ldots, k-1\}$, indexed by $b$ \\
        {$E$} & Set of exams, indexed by $e$\\ 
        \underline{Parameters}  \\
        $k$ & Desired number of exam blocks \\
        $c_{ij}$ & The number of students co-enrolled in both exam $i$ and exam $j$ \\
        \underline{Decision Variables}  \\
        $x_{eb}$ & = 1 if exam $e$ is assigned to block $b$; 0 otherwise\\
        $y_{ij}$ & = 1 if exam $i$ and exam $j$ are assigned to the same block; 0 otherwise\\
        $z_{ij}$ &= 1 if exam $i$ and exam $j$ are assigned to neighboring blocks; 0 otherwise\\
        \bottomrule
    \end{tabular}
\end{table*}

\subsection{Min-Conflict IP Formulation} \label{sec:min-ip}

The \textsc{Min\_Conflict\_IP} problem is presented as follows.

\begin{align}
    \min~~~ & \sum\limits_{i, j \in E} c_{ij} \ y_{ij} \label{ass-min-obj} \\
    \text{s.t.}~~~ & \sum_{b \in B}x_{eb} = 1, \forall ~ e \in E \label{ass-min-c1} \\ 
    & x_{ib} + x_{jb} \leq 1+y_{ij}, \forall  ~i, j \in E, ~ b \in B \label{ass-min-c2}\\
    & x_{eb} \in \{0, 1\}, \forall ~ e \in E, b \in B \label{ass-min-c3}\\
    & y_{ij} \in \{0, 1\}, \forall ~ i, j \in E\label{ass-min-c4}
\end{align}

Objective function (\ref{ass-min-obj}) ensures that the total number of conflicts is minimized. Constraints (\ref{ass-min-c1}) ensure that each exam is assigned to exactly one block. Constraints (\ref{ass-min-c2}) ensure that $y_{ij} = 1$ when both exam $i$ and exam $j$ are assigned to the same block and $=0$ otherwise. Constraints (\ref{ass-min-c3}) and (\ref{ass-min-c4}) are binary constraints for the variables.

\subsection{Zero-Conflict IP Formulation} \label{sec:zero-ip}

The \textsc{Zero\_Conflict\_IP} problem is formulated similar to the \textsc{Min\_Conflict\_IP} problem but with an updated TSP objective function and additional constraints.

The objective function is updated to:
\begin{align}
    \min~~~ & \sum\limits_{i, j \in E} c_{ij} \ z_{ij} \label{ass-zero-obj}
\end{align}

And the following constraints are added:
\begin{align}
    \text{s.t.}~~~ & \sum\limits_{i, j \in E} c_{ij} \ y_{ij} = 0 \label{ass-zero-c1} \\
    & x_{ib} + x_{jb+1} \leq 1 + z_{ij}, \forall ~i,j \in E, ~b \in \{0, \ldots, k-2\} \label{ass-zero-c2}\\
    & z_{ij} \in \{0, 1\}, \forall ~ i, j \in E \label{ass-zero-c3}
\end{align}

Objective function (\ref{ass-zero-obj}) minimizes the total number of conflicts of neighbouring blocks. Constraint (\ref{ass-zero-c1}) ensures that there is no conflict within a block. Constraints (\ref{ass-zero-c2}) and (\ref{ass-zero-c3}) define the $z_{ij}$ variables and impose the binary constraints, respectively.

\section{The Max-Clique Problem}\label{sec:max_clique_ip}
\begin{table*}[!hp]\scriptsize
    \centering
    \caption{Notations for the Max-Clique Problem}
    \begin{tabular}{l l }
        \toprule
        \textbf{Notation} & \textbf{Description} \\
        \midrule
        \underline{Index Sets}  \\
        {$\mathrm{NN}_e$} & Set of non-neighbors of exam $e$ \\
        \underline{Parameters}  \\
        $h_e$ & The number of non-neighbors of exam $e$ \\
        \underline{Decision Variables}  \\
        $x_e$ & = 1 if exam $e$ is included in the clique; 0 otherwise \\
        \bottomrule
    \end{tabular}
\end{table*}

\subsection{The Max-Clique IP Formulation}

The \textsc{Max\_Clique\_IP} problem is presented as follows.

\begin{align}
    \max~~~ & \sum\limits_{e \in E} x_e \label{cli-max-obj} \\
    \text{s.t.}~~~ & \revision{\sum\limits_{e' \in \mathrm{NN}_e} x_{e'} \leq h_e (1-x_e)}, \forall e \in E \label{cli-max-c1}\\
    & \revision{x_e \in \{0, 1\}}, \forall e \in E \label{cli-max-c2}
\end{align}

Objective function (\ref{cli-max-obj}) minimizes the number of exams included in the clique. Constraints (\ref{cli-max-c1}) ensure that if an exam is included in the clique, then none of its non-neighbors can be included. Constraints (\ref{cli-max-c2}) enforce binary values for the decision variables.

\section{The Block Sequencing Problem} \label{sec:sequencing_ip}

The \textsc{Block\_Sequencing\_IP} problem is presented as follows.

\begin{table*}[!ht]\scriptsize
    \centering
    \caption{Notations for the Block Sequencing Problem}
    \label{tab:param2}
    \begin{tabular}{p{3cm}p{12cm}p{1cm}p{6cm}}
        \toprule
        \textbf{Notation} & \textbf{Description} \\
        \midrule
        \underline{Sets}  \\
        {$B$} & Set of exam blocks, $B := \{0, \ldots, |S|-1\}$ \\
        {$S$} & Set of time slots available, $S= \{0, \ldots, |S|-1\}$\\ 
        {$Tripleday$} & Subset of time slots that index the start of triples on the same day, $Tripleday \subset S$  \\ 
        {$Triple24hr$} & Subset of time slots that index the start of triples in 24 hours, $Triple24hr \subset S$, \revision{with $Tripleday ~\cap ~ Triple24hr = \emptyset$} \\ 
        {$Tripleslots$} & {$Tripleday ~ \cup ~ Triple24hr$}\\
        {$B2B_{evemorn}$} & Subset of time slots that index the start of evening-morning back-to-backs, $B2B_{evemorn} \subset S$ \\ 
        {$B2B_{other}$} & Subset of time slots that index the start of the other types of back-to-backs, $B2B_{other} \subset S$, \revision{with $B2B_{evemorn} ~\cap~ B2B_{other} = \emptyset$}\\
        $Earlyslots$ & Subset of slots earlier than a slot cutoff, $Earlyslots \subseteq S$ \\
        $LargeB$  & Subset of blocks that have enrollment larger than a size cutoff, $LargeB \subseteq B$ \\
        \underline{Parameters}  \\
        {$n_s$} & Next time slot of time slot $s$ in the cycle, $n_s:=(s+1) \mod |S|$  \\ 
        {$\beta$} & Penalty weight for three exams in 24 hours \\ 
        {$\alpha$} & Penalty weight for three exams on the same day (Having one exam in the afternoon, one in the evening, and one on the next morning would count as having three exams in 24 hours, not three exams on the same day)  \\ 
        {$\gamma_1$} & Penalty weight for evening-morning back-to-back exams \\
        {$\gamma_2$} & Penalty weight for the other types of back-to-back exams \\
        $\delta$ & Penalty weight for three exams in four slots\\
        $p_{ij}$ & The number of students co-enrolled in both block $i$ and block $j$ \\ 
        {$t_{ij\ell}$} & The number of students co-enrolled in all blocks $i,j$ and $\ell$ \\
        \underline{Decision Variables}  \\
        {$x_{ij\ell s}$} & = 1 if block $i$ is placed in slot $s$, block $j$ is placed in slot $(s+1) \mod |S|$, and block $\ell$ is placed in slot $(s+2) \mod |S|$; 0 otherwise\\ 
        {$y_{ij\ell}$} & = 1 if block $i, j, \ell$ appear in a sequence where each pair of consecutive blocks is directly back-to-back; 0 otherwise \\
        {$z_{ij\ell m}$} & = 1 if block $i, j, \ell, m$ appear in a sequence where each pair of consecutive blocks is directly back-to-back; 0 otherwise\\
        \bottomrule
    \end{tabular}
\end{table*}

\subsection{Objective Function}
\begin{align}\label{seq-obj} 
    \min~~~ &\sum\limits_{i, j, \ell \in B, s \in B2B_{evemorn}} \gamma_1 \ p_{ij} \ x_{ij\ell s} \\ \nonumber
    & + \sum\limits_{i, j, \ell \in B, s \in B2B_{other}} \gamma_2 \  p_{ij} \ x_{ij\ell s} \\ \nonumber
    &
    + \sum\limits_{i, j, \ell \in B, s \in Tripleday} \alpha \ t_{ij\ell} \ x_{ij\ell s}  \\ \nonumber
    &
    + \sum\limits_{i, j, \ell \in B, s \in Triple24hr} \beta \ t_{ij\ell} \ x_{ij\ell s} \\ \nonumber
    & +\sum\limits_{i, j, \ell, m \in B} \delta \ (t_{ij\ell}+t_{i\ell m}) \ z_{ij\ell m}
\end{align}

Objective \revision{function \eqref{seq-obj}} minimizes the weighted sum of different higher-level conflicts, including back-to-backs, three exams in a row, and three exams in four slots.

\subsection{Constraints}
\begin{align}
    \text{s.t.}~~~ & \sum\limits_{j,\ell \in B, s \in S}x_{ij\ell s}  = 1, \forall i \in B \label{seq-c1} \\ 
    & \sum\limits_{i,\ell \in B, s \in S}x_{ij\ell s}  = 1, \forall j \in B \label{seq-c2}\\
    & \sum\limits_{i,j \in B, s \in S}x_{ij\ell s}  = 1, \forall \ell \in B \label{seq-c3}\\
    & \sum\limits_{i,j,\ell \in B}x_{ij\ell s}  = 1, \forall s \in S \label{seq-c4}\\
    & x_{ii\ell s} = 0, \forall i,\ell \in B, s \in S \label{seq-c5}\\
    & x_{ijis} = 0, \forall i,j \in B, s \in S \label{seq-c6}\\
    & x_{ijjs} = 0, \forall i,j \in B, s \in S \label{seq-c7} \\
    & \sum\limits_{i \in B}x_{ij\ell s}  = \sum\limits_{m \in B}x_{j \ell m n_s}, \forall j,\ell \in B, s \in S \label{seq-c8}\\
    & y_{ij\ell} = \sum\limits_{s \in Tripleslots}x_{ij\ell s} , \forall i,j,\ell \in B \label{seq-c9}\\
    & z_{ij\ell m} \geq y_{ij\ell} + y_{j\ell m} - 1, \forall i,j,\ell,m \in B \label{seq-c10}\\
    & x_{ij\ell s}  \in \{0, 1\}, \forall i, j, \ell \in B, s \in S \label{seq-c11}\\
    & y_{ij\ell} \in\{0, 1\}, \forall i, j, \ell \in B \label{seq-c12}\\
    & z_{ij\ell m} \in\{0, 1\}, \forall i, j, \ell, m \in B \label{seq-c13}
\end{align}

Constraints (\ref{seq-c1}) - (\ref{seq-c3}) ensure that each block appears as the first, second, and third of a triple exactly once. Constraint (\ref{seq-c4}) ensures that each time slot is assigned to exactly one block. Constraints (\ref{seq-c5}) - (\ref{seq-c7}) prevent repetition in group sequences. Constraints (\ref{seq-c8}) define the triplet time slot continuity (e.g. if $x_{1,2,3,1} = 1$, then $x_{2,3,\ell,2} = 1$ for some block $\ell$). Constraints (\ref{seq-c9})-(\ref{seq-c13}) define the $y_{ij\ell}$ and $z_{ij\ell m}$ variables and impose the binary constraints.

Front-loading constraints (\ref{seq-c14}) can be added to keep large exams out of the last couple of slots:
\begin{align}
    \sum\limits_{j,\ell \in B, s \in Earlyslots}x_{ij\ell s}  = 1, \forall ~ i \in LargeB \label{seq-c14}
\end{align}

\section{The Post-Processing Algorithm}
\label{sec:post_processing_appendix}

Algorithm \ref{alg:local_search} details the MIP-based local search algorithm used in the Post-Processing Stage. 

\begin{table*}[!ht]\scriptsize
    \centering
    \caption{Notations for the Schedule\_IP Problem}
    \label{tab:param3}
    \begin{tabular}{p{3cm}p{12cm}p{1cm}p{6cm}}
        \toprule
        \textbf{Notation} & \textbf{Description} \\
        \midrule
        \underline{Sets}  \\
        {$E'$} & A subset of exams, $E' \subseteq E$\\ 
        {$S'$} & A subset of available time slots, $S' \subseteq S$\\ 
        $B2B$ & Set of time slots that index the start of back-to-backs, $B2B \subset S'$\\
        $2i3$  & Set of time slots that index the start of two-in-threes, $2i3 \subset S'$\\
        \underline{Parameters}  \\
        $\lambda_1$ & Relative weight between conflicts and back-to-backs\revision{; default is 1000}\\
        $\lambda_2$ & Relative weight between two-in-threes and back-to-backs\revision{; default is 0.5}\\
        $\zeta_s$ & Front-loading penalty for large exams scheduled to slot $s$\revision{; default is $s^{300}$}\\
        $other\_conf_{es}$ & Total number of conflicts between exam $e$ and any exam $e' \in E \setminus \{e\}$ scheduled to slot s (i.e., $x_{e's} = 1$)\\
        $other\_B2B_{es}$ &  Total number of back-to-backs between exam $e$ and any exam $e' \in E \setminus \{e\}$ scheduled to slot s \\
        $other\_2i3_{es}$ & Total number of two-in-threes between exam $e$ and any exam $e' \in E \setminus \{e\}$ scheduled to slot s \\
        \underline{Decision Variables}  \\
        $x_{es}$ & = 1 if exam $e \in E'$ is in slot $s$; 0 otherwise\\
        $y_{ij}$ & = 1 if exams $i, j \in E'$ are in the same slot; 0 otherwise\\
        $z_{ij}$ & = 1 if exams $i, j \in E'$ are in consecutive slots; 0 otherwise\\
        $w_{ij}$  & = 1 if exams $i, j \in E'$ are two slots apart; 0 otherwise\\
        \bottomrule
    \end{tabular}
\end{table*}

The \textsc{Schedule\_IP} problem used in both the Post-Processing local search algorithm and the Layer-Cake heuristic algorithm is presented below.

\subsection{Objective Function}
\begin{align}\label{post-obj} 
    \min~~~ & \lambda_1\ (\sum\limits_{i, j \in E'}0.5 \ c_{ij} \ y_{ij} + \sum\limits_{e \in E', s \in S'} other\_conf_{es} \ x_{es}) + \\\nonumber
    & \sum\limits_{i, j \in E'}c_{ij}\ z_{ij} + \sum\limits_{e \in E', s \in S'} other\_B2B_{es}\ x_{es} + \\\nonumber
    & \lambda_2\ (\sum\limits_{i, j \in E'}c_{ij} \ w_{ij} + \sum\limits_{e \in E', s \in S'} other\_2i3_{es} \ x_{es})+ \\\nonumber
    & \sum\limits_{e \in E', s \in S' \setminus Earlyslots}\zeta_{s} \ x_{es}
\end{align}

Objective function (\ref{post-obj}) minimizes the weighted sum of conflicts, back-to-backs, and two-in-threes.  

\subsection{Constraints}
\begin{align}
\text{s.t.}~~~     & \sum\limits_{s \in S'}x_{es} = 1, \forall e \in E' \label{post-c4}
    \\
    & y_{ij} \geq x_{is} + x_{js} - 1, \forall i,j \in E', s \in S' \label{post-c1}
    \\ 
    & z_{ij} \geq x_{is} + x_{js+1} - 1, \forall i,j \in E', s \in B2B \label{post-c2}
    \\ 
    & w_{ij} \geq x_{is} + x_{js+2} - 1, \forall i,j \in E', s \in 2i3 \label{post-c3}
    \\ 
    & x_{es} \in\{0, 1\}, \forall e \in E', s \in S' \label{post-c5}
    \\
    & y_{ij} \in\{0, 1\}, \forall i, j \in E' \label{post-c6}
    \\
    & z_{ij} \in\{0, 1\}, \forall i, j \in E' \label{post-c7}
    \\
    & w_{ij} \in\{0, 1\}, \forall i, j \in E' \label{post-c8}
\end{align}

Constraints (\ref{post-c4}) ensures that every exam is assigned to exactly one slot. Constraints (\ref{post-c1})-(\ref{post-c8}) define the $y_{ij}, z_{ij}$, and $w_{ij}$ variables and impose the binary constraints.

\begin{algorithm}[!ht] \scriptsize
 \caption{A MIP-based Local Search Algorithm for Post-Processing}\label{alg:local_search}
\begin{algorithmic}[1]
\State \textbf{Input:} $base\_sched, S$ \Comment{A given complete schedule that maps the set of exams to the set of time slots $S$}
\State $n \gets 25$ \Comment{Number of exams to be rescheduled}
\State $f_{step},b_{step}  \gets \textsc{Round}(0.4n), \textsc{Round}(0.2n)$
\State $sched \gets base\_sched$
    \State $i \gets 0$ \Comment{Starting index of exams to be rescheduled}
    \While{$i < |sched|$}
        \State $best\_score \gets \textsc{Get\_Score}(sched)$ * 
        \State $E \gets \textsc{Sort\_By\_Density}(sched)$\Comment{Sort the schedule by the density of bad events}
        \State ${new\_sched} \gets \textsc{Reschedule}(sched, E[i:i+n])$  \Comment{Reschedule the next $n$ exams from index $i$}
        \State $new\_score \gets \textsc{Get\_Score}(new\_sched)$
        \If{$best\_score \leq new\_score$} 
            \State $i \gets i + f_{step}$ \Comment{Increase the starting index by a larger step if no improvement was found}
                \Else 
            \State $i \gets max(0, i - b_{step})$ \Comment{Decrease the starting index by a smaller step if a better schedule was found}
            \State $sched \gets new\_sched$
        \EndIf
    \EndWhile

\State \textbf{Output:} $sched$

\vskip6pt
\Procedure{Get\_Score\_Individual}{$e, s$}
    \State \textbf{return} $\lambda_1 \ other\_conf_{es} + other\_{B2B}_{es} + \lambda_2 \ other\_2i3_{es}$
\EndProcedure
\Statex

\Procedure{Get\_Score}{$sched$} \Comment{Given a schedule, i.e., a mapping from exams to time slots}
    \State \textbf{return} $\sum_{(e, s) \in sched} \Call{Get\_Score\_Individual}{e, s}$
\EndProcedure
\Statex

\Procedure{Reschedule}{$sched, E'$} \Comment{Given a complete schedule and a set of exams to reschedule}
    \State Solve \textsc{Schedule\_IP}$(E' = E', S' = S)$ and enforce the following additional constraints:
    \State \quad 
    $x_{es} = 1, \forall (e, s) \in sched: e \notin E'$
    \State \textbf{return} $\{(e,s)\mid x_{es}=1\}$
    
\EndProcedure
\Statex
\end{algorithmic}

\end{algorithm}

\section{The Layer-Cake Heuristic Algorithm}\label{sec:layer_cake_appendix}

Algorithm \ref{alg:layer_cake} details the MIP-based Layer-Cake heuristic algorithm.

To speed up the algorithm, we construct a warm-start solution for \textsc{Schedule\_IP} at each iteration of Layer-Cake as follows:
\begin{enumerate}
    \item Exams that are re-entering are assigned time slots based on the results from the previous layer.
    \item New exams, not present in the last layer, are first ordered by size and then assigned time slots greedily. For each exam, the time slot that results in the lowest score, as determined by \textsc{Get\_Score\_Individual}, is selected, considering the exams that have already been placed.
\end{enumerate}

\begin{algorithm}[!ht] \scriptsize
\caption{A MIP-based Layer-Cake Heuristic Algorithm}\label{alg:layer_cake}
\begin{algorithmic}[1]
\State \textbf{Input:}  $S_2$ \Comment{Set of available time slots}
\State $n_2 \gets 15000$ \Comment{The desired number of student-exam pairs per layer}
\State $sched, layer \gets \{\}, \emptyset$ \Comment{Initializing the schedule}
\State $n_3 \gets \textsc{Round}(0.3n_2)$  \Comment{Desired number of student-exam pairs that enter the layer for a second time}
\State $E \gets \textsc{Sort\_By\_Size}(E)$ \Comment{Set of exams ordered by size descending}
\State $layer \gets \emptyset$
\For{$i = n_2;\ i>0;\ i = i- \textsc{Size}(E[0])$}
\Comment{Construct the first layer}
    \State $layer \gets layer \cup E[0]$
    \State $E \gets E \setminus E[0]$
  \EndFor
\State $sched \gets \textsc{Layer\_Schedule}(sched, \emptyset, layer, S_2)$

\While{$E \neq \emptyset$}
    \State $old\_layer, overlap \gets layer, \emptyset$
    \For{$i = n_3;\ i>0;\ i = i- \textsc{Size}(old\_layer[-1])$}  \Comment{Append set of exams to overlap from the previous layer}
            \State $overlap \gets overlap \cup old\_layer[-1]$
            \State $old\_layer \gets old\_layer \setminus old\_layer[-1]$
      \EndFor
    \State $layer \gets overlap$
    \For{$i = n_2-n_3;\ i>0;\ i = i- \textsc{Size}(E[0])$}
    \Comment{Append set of exams to the current layer from the remaining exams}
        \State $layer \gets layer \cup E[0]$
        \State $E \gets E \setminus E[0]$
      \EndFor
    \State $sched \gets \textsc{Layer\_Schedule}(sched, overlap, layer, S_2)$
\EndWhile
\State \textbf{Output:} $sched$
\vskip6pt

  \revision{\Procedure{Layer\_Schedule}{$sched,\,overlap,\,E',\,S'$}
  \Comment{Given an (incomplete) schedule, a set of re‐entering exams, a set of exams to schedule, and a set of available slots}
\State {\scriptsize\textbf{--- Build warm‐start solution }$x^0$}
  \ForAll{$(e,s)\in sched$ \textbf{with} $e\in overlap$}
    \State set $x^0_{es} \gets 1$  
  \EndFor
  \State $E_{\mathrm{new}} \gets E' \setminus overlap$
  \State sort $E_{\mathrm{new}}$ by descending size
  \ForAll{$e\in E_{\mathrm{new}}$}
    \State pick 
      $s^* \;=\;\arg\min_{s\in S'}\textsc{Get\_Score\_Individual}(e,s)$
    \State set $x^0_{e\,s^*}\gets 1$
  \EndFor
  \State {\scriptsize\textbf{--- MIP solve with warm start}}
  \State Solve \textsc{Schedule\_IP}$(E'=E',S'=S')$ providing $x^0$ as the initial solution and enforce the following additional constraints: 
    \State \quad
    $x_{es} = 1, \forall (e, s) \in sched: e \notin overlap$
  \State \Return $\{(e,s)\mid x_{es}=1\}$
\EndProcedure}
\Statex
\end{algorithmic}
\end{algorithm}

\section{Experimental Setting} \label{sec:setting_appendix}

By default, we solve the \textsc{Min\_Conflict\_IP} problem for Block Assignment. For Block Sequencing, Post-Processing, and Layer-Cake, the following parameters are used unless otherwise stated. \revision{The ratios of Conflicts to B2Bs, Triples to B2Bs, 3-in-4s to B2Bs, and 2-in-24hrs to B2Bs are 1000, 10, 5, and 0.5, respectively.} For front-loading parameters, the default cutoff for early time slots is 23, and the large exam size cutoff is 300. No time slot is excluded from the available set by default. In the experiments presented, we set a time limit of 1500 seconds for the \textsc{Block\_Sequencing\_IP} model, 600 seconds for the Post-Processing Algorithm, and 1500 seconds for each \textsc{Schedule\_IP} in Layer-Cake. In practice, we would extend the time limits reasonably large to allow the models to solve to optimality.

\section{Nottingham Dataset Benchmark} \label{sec:Nottingham_appendix}

In order to compare our framework with a selection of existing exam timetabling algorithms, we evaluate an adapted version of our model using the University of Nottingham dataset introduced in \cite{burke1996memetic}. Specifically, the Nottingham \emph{a} variant aims to minimize instances of students having back-to-back exams on the same day, whereas the Nottingham \emph{b} variant additionally seeks to minimize a weighted total of within-day and overnight back-to-backs. We selected this dataset not only because of the similarity in the objective functions, but also to demonstrate our framework's capability to handle capacitated examination timetabling problems. Additional constraints are detailed in \cite{qu2009survey}.

In adapting the Zero-GtSP model for the Nottingham dataset, the following modifications and constraints were introduced:
\begin{itemize}

    \item \textbf{Coincident Courses Grouping:} Courses that share identical student enrollments or are intended to be scheduled concurrently are aggregated into a single exam within the model.

    \item \textbf{Default Block Ordering as Time Slots:} To enforce the specific time slot and sequential constraints in all stages, we assume that the default ordering of blocks generated during the block assignment stage directly corresponds to the sequence of time slots. Consequently, the time slot set \(S\) is used interchangeably with the block set \(B\) in the subsequent constraints.

    \item \textbf{Capacity Constraints:} Room and seating capacities are strictly enforced for each exam slot, ensuring that no block exceeds its predefined limits. Let $C = 1550$ denote the capacity limit per time slot and let $q_e$ represent the size of exam $e$. The capacity constraint is then expressed as
    \[
    \sum_{e \in E} q_e\, x_{es} \leq C, \quad \forall\, s \in S.
    \]

    \item \textbf{Specific Time Slot Constraints:} Certain exams are restricted to specific time slots in accordance with institutional requirements (e.g., some exams must be scheduled in the morning or afternoon). Let $S_e \subset S$ be the set of permissible time slots for exam $e$. The corresponding constraint is given by
    \begin{align}
         \sum_{s \in S_e} x_{es}  = 1, \quad \forall\, e \in E. \label{eq:allowed-block-constraint}
    \end{align}

    \item \textbf{Sequential Dependencies:} Exams that require a predefined order (i.e., one exam must precede another) are handled via sequencing constraints. Let $\mathcal{P} = \{(e_1, e'_1), (e_2, e'_2), \ldots\}$ denote the set of exam pairs with precedence requirements. Define
    \begin{align}
         y_e &= \sum_{s \in S} s\ x_{es}, \quad \forall\, e \in E,
    \end{align}
    and impose the ordering condition:
    \begin{align}
         y_e + 1 &\leq y_{e'}, \quad \forall\, (e, e') \in \mathcal{P}.
    \end{align}

    \item \textbf{Non-Overlapping Constraints:} Exams identified as incompatible (e.g., due to shared instructors, equipment, or other resource limitations) are scheduled in non-overlapping time slots. Let $\mathcal{G} \subset E$ denote a group of such exams that must be assigned distinct time slots. The corresponding constraint is:
    \[
    \sum_{e \in \mathcal{G}} x_{es} \leq 1, \quad \forall\, s \in S.
    \]
\end{itemize}

According to \cite{burke1996memetic}, only the zero-conflict constraint and the slot capacity constraint are treated as hard constraints, while the remaining constraints described above may be considered soft constraints.

\revision{
For this set of experiments, we used a time limit of 10 hours for block assignment, 3600 seconds for the \textsc{Block\_Sequencing\_IP} model, and 1800 seconds for the Post-Processing Algorithm.
}

Table~\ref{tab:nottingham} compares the solution outcomes of established exam timetabling algorithms \citep{qu2009survey} with those obtained using our adapted Zero-GtSP model. Figure~\ref{fig:nottingham} illustrates the distribution of exam schedules generated by our approach. The strong performance of our method relative to the literature underscores the flexibility and robustness of our framework, even under stringent capacity constraints. Moreover, our framework is capable of addressing higher-order conflicts that other algorithms are unable to manage.

\begin{table}[ht]
\centering
\caption{Comparison of Objective Value for Nottingham Dataset \emph{a} and \emph{b}}
\scriptsize
\begin{adjustbox}{width=\textwidth}
\begin{tabular}{lcccccc}
\toprule
Nottingham & \makecell{\cite{burke1996memetic}} & \makecell{\cite{di2001tabu}} & \makecell{\cite{caramia2000new}} & \makecell{\cite{merlot2003hybrid}} & \makecell{\cite{abdullah2007tabu}} & Adapted Zero-GtSP (Ours) \\
\midrule
\emph{a} (26 slots) & 53  & 11  & 44  & 2   & 18  & 17 \\
\emph{a} (23 slots) & 269 & 123 & --  & 88  & --  & 0  \\
\midrule
Nottingham & \makecell{\cite{burke1999multistage}} & \makecell{\cite{di2001tabu}} & \makecell{\cite{merlot2003hybrid}} & \makecell{\cite{burke2004solving}} & \makecell{\cite{burke2004time}} & Adapted Zero-GtSP (Ours)\\
\midrule
\emph{b} (26 slots) & --  & --  & --  & --  & --  & 244\\
\emph{b} (23 slots) & 519 & 751 & 401 & 545 & 384 & 330 \\
\bottomrule
\end{tabular}
\end{adjustbox}
\label{tab:nottingham}
\end{table}

\begin{figure}[htbp!]
    \centering
    \begin{subfigure}[b]{0.45\textwidth}
        \centering
        \includegraphics[width=\textwidth]{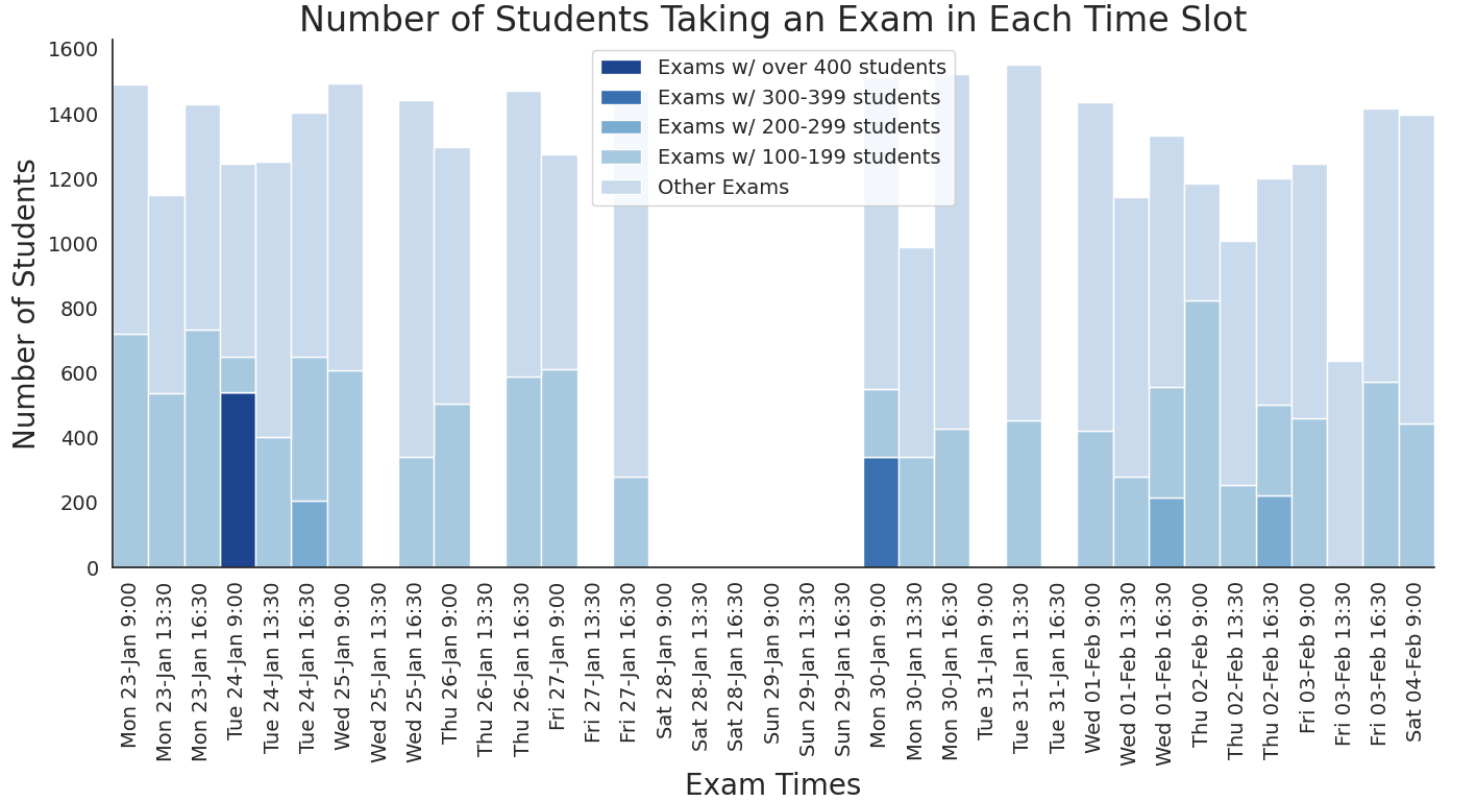}
        \caption{Nottingham \emph{a} with 26 slots}
    \end{subfigure}
    \hspace{1em}
    \begin{subfigure}[b]{0.45\textwidth}
        \centering
        \includegraphics[width=\textwidth]{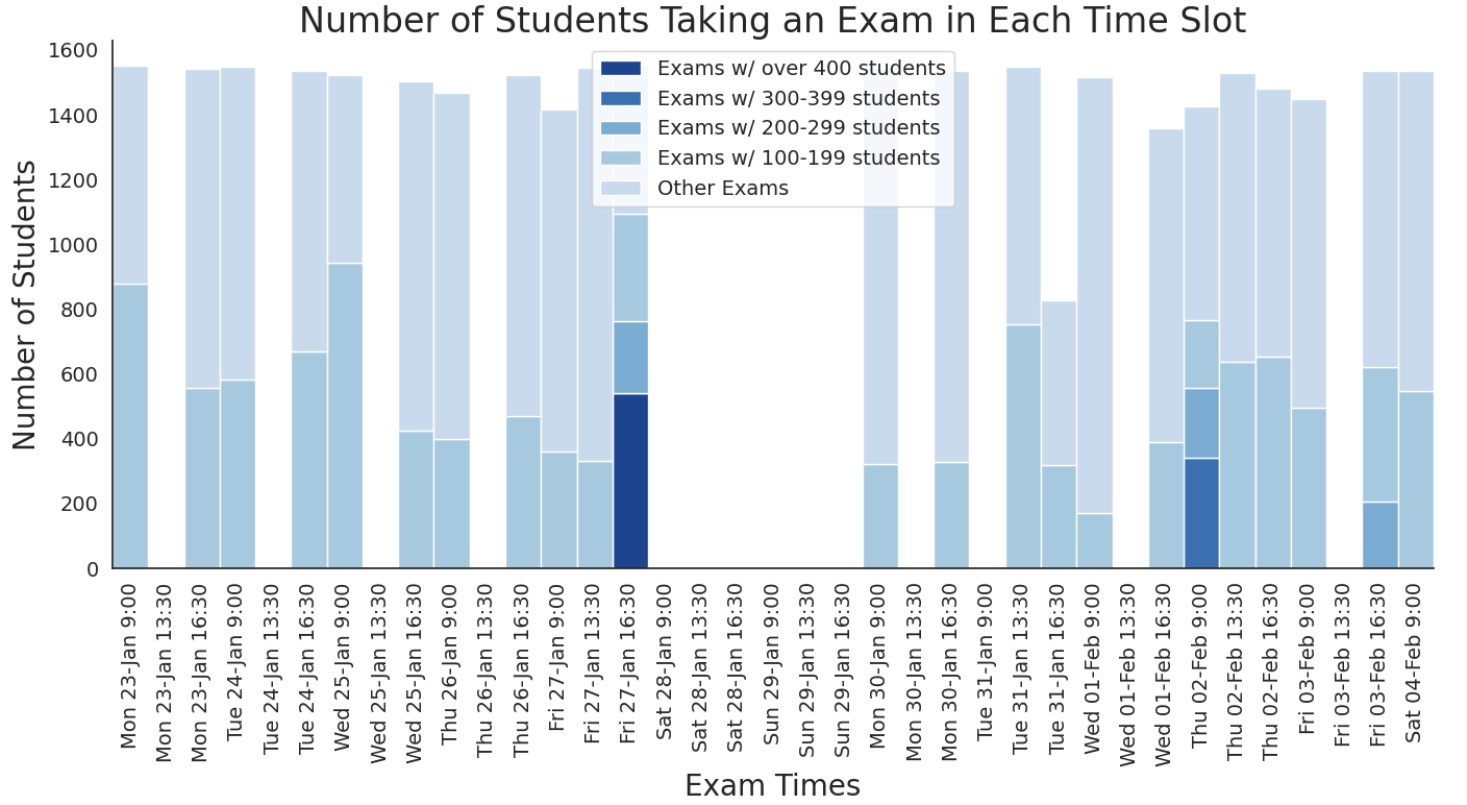}
        \caption{Nottingham \emph{a} with 23 slots}
    \end{subfigure}

    \vspace{1em} 

    \begin{subfigure}[b]{0.45\textwidth}
        \centering
        \includegraphics[width=\textwidth]{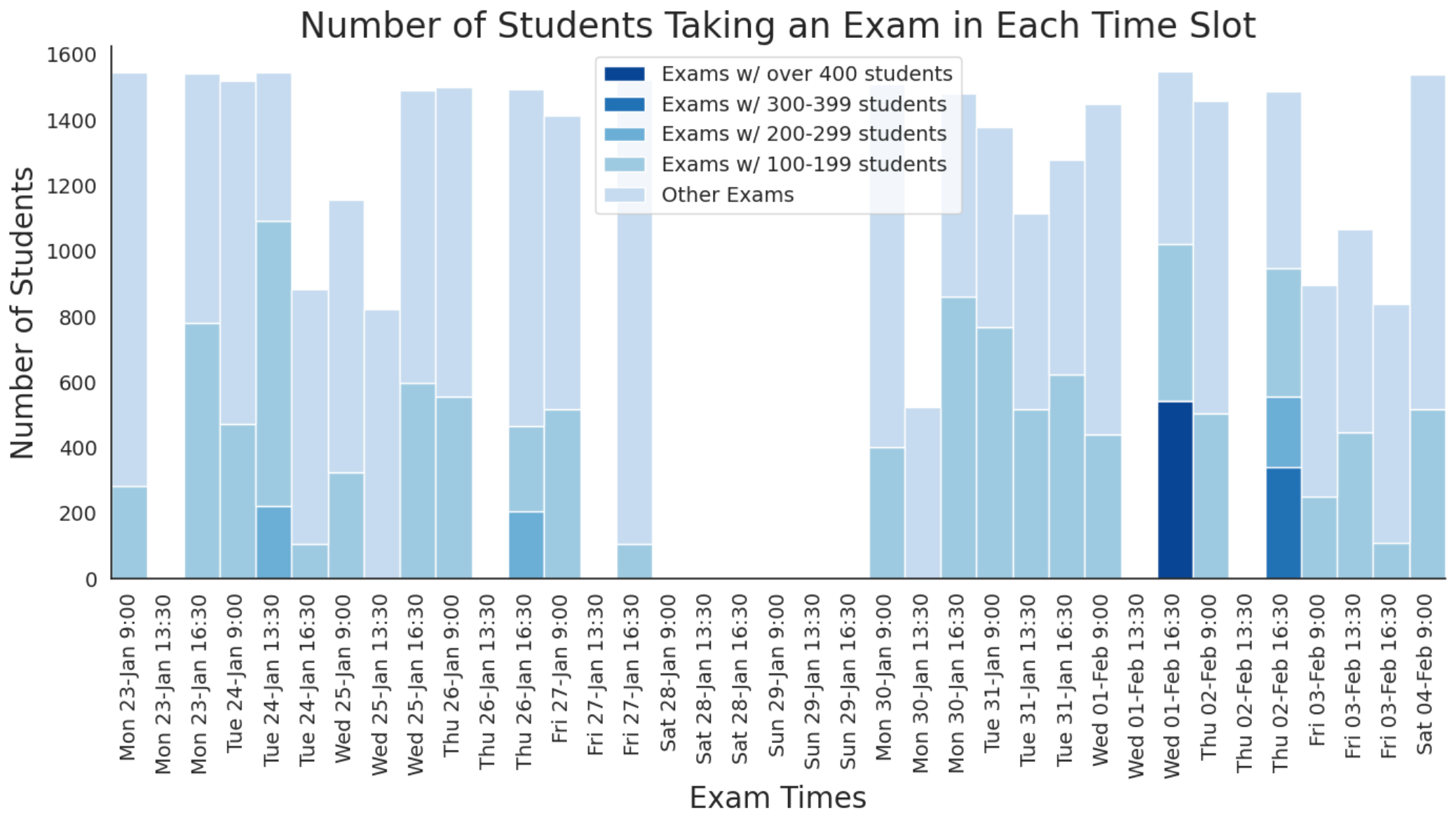}
        \caption{Nottingham \emph{b} with 26 slots}
    \end{subfigure}
    \hspace{1em}
    \begin{subfigure}[b]{0.45\textwidth}
        \centering
        \includegraphics[width=\textwidth]{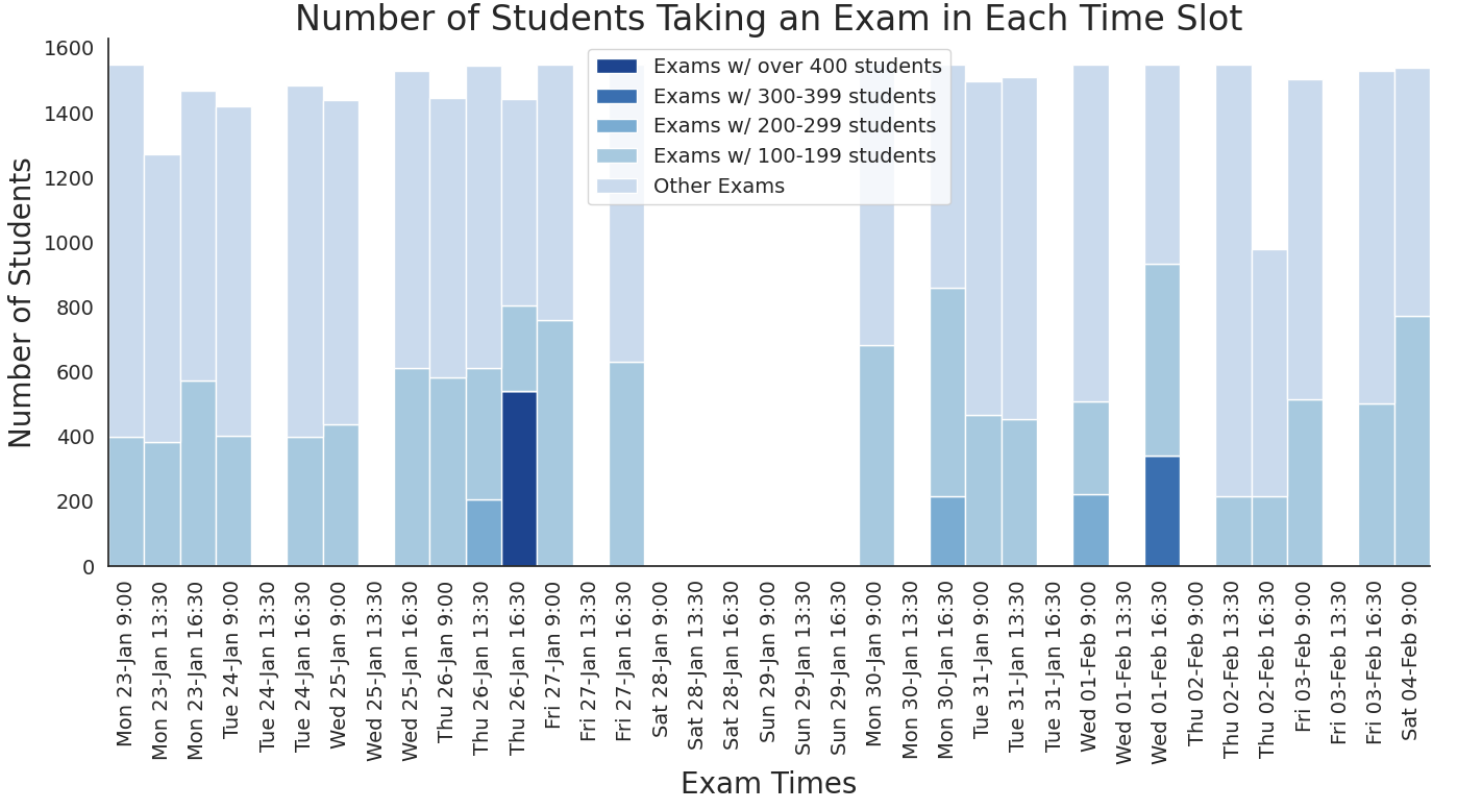}
        \caption{Nottingham \emph{b} with 23 slots}
    \end{subfigure}

    \caption{Distribution of Exams Scheduled by Adapted Zero-GtSP on Nottingham Dataset \emph{a} and \emph{b}}
    \label{fig:nottingham}
\end{figure}

\end{APPENDICES}

\end{document}